\documentclass[11pt,
]{article}

\usepackage[all]{xy}

\usepackage[latin1]{inputenc}
\usepackage{amsfonts}
\usepackage{amsmath}
\usepackage{amssymb}
\usepackage{amscd}
\usepackage{amsthm}
\usepackage{indentfirst}
\usepackage[hmargin=2.8cm
,vmargin=3.9cm
]{geometry}

\usepackage{faktor}

\usepackage{epsfig}

\newtheorem{thm}{Theorem}[section]
\newtheorem{lemma}[thm]{Lemma}
\newtheorem{conj}[thm]{Conjecture}
\newtheorem{prop}[thm]{Proposition}
\newtheorem{coroll}[thm]{Corollary}

\theoremstyle{definition}

\newtheorem{defin}[thm]{Definition}
\newtheorem{rem}[thm]{Remark}
\newtheorem{exam}[thm]{Example}

\newtheorem{notation}[thm]{Notation}

\newcommand{\R}{{\mathbb{R}}}

\newcommand{\T}{{\mathbb{T}}}
\newcommand{\Z}{{\mathbb{Z}}}
\newcommand{\N}{{\mathbb{N}}}
\newcommand{\C}{{\mathbb{C}}}

\newcommand{\cA}{{\mathcal{A}}}
\newcommand{\cB}{{\mathcal{B}}}
\newcommand{\cC}{{\mathcal{C}}}

\newcommand{\cF}{{\mathcal{F}}}

\newcommand{\cO}{{\mathcal{O}}}
\newcommand{\cP}{{\mathcal{P}}}

\newcommand{\cS}{{\mathcal{S}}}
\newcommand{\cT}{{\mathcal{T}}}
\newcommand{\cU}{{\mathcal{U}}}
\newcommand{\cV}{{\mathcal{V}}}

\newcommand{\fc}{{:\ }}

\newcommand{\ol}{\overline}
\newcommand{\wt}{\widetilde}
\newcommand{\wh}{\widehat}

\newcommand{\tb}{\textbf}

\DeclareMathOperator{\spn}{span}

\DeclareMathOperator{\im}{im}
\DeclareMathOperator{\id}{id}

\DeclareMathOperator{\pr}{pr}

\DeclareMathOperator{\Lip}{Lip}

\DeclareMathOperator{\const}{const}

\DeclareMathOperator{\std}{std}

\begin{document}

\title{Constraints on symplectic quasi-states}

\author{Adi Dickstein\footnote{Partially supported by the Milner Foundation and by the Israel Science Foundation grant 1102/20.}\ \ and Frol Zapolsky}

\date{}

\setcounter{tocdepth}{2}

\renewcommand{\labelenumi}{(\roman{enumi})}

\maketitle

\begin{abstract}
  We prove that given a closed connected symplectic manifold equipped with a Borel probability measure, an arbitrarily large portion of the measure can be covered by a symplectically embedded polydisk, generalizing a result of Schlenk. We apply this to constraints on symplectic quasi-states. Quasi-states are a certain class of not necessarily linear functionals on the algebra of continuous functions of a compact space. When the space is a symplectic manifold, a more restrictive subclass of symplectic quasi-states was introduced by Entov--Polterovich. We use our embedding result to prove that a certain `soft' construction of quasi-states, which is due to Aarnes, cannot yield nonlinear symplectic quasi-states in dimension at least four.
\end{abstract}

\tableofcontents

\section{Introduction and results}\label{s:intro}

\subsection{Quasi-states and the main result}\label{ss:QS_main_result}

A \tb{quasi-state} on a compact Hausdorff space $X$ is a functional $\zeta \fc C(X) \to \R$, where $C(X)$ is the algebra of continuous real-valued functions on $X$, which satisfies:
\begin{itemize}
  \item \tb{(normalization):} $\zeta(1) = 1$;
  \item \tb{(monotonicity):} $\zeta(f) \leq \zeta(g)$ whenever $f \leq g$;
  \item \tb{(quasi-linearity):} for each $f$, $\zeta$ is linear on the corresponding singly-generated subalgebra $C(f) = \{h\circ f\,|\,h\fc \R \to \R\text{ continuous}\}$.
\end{itemize}
\begin{rem}\label{rem:Riesz_rep}
\begin{enumerate}
  \item By the Riesz representation theorem, $\mu \mapsto \int_X\cdot\,d\mu$ is a bijection between the collections of regular\footnote{That is $\mu(E) = \sup\{\mu(K)\,|\,K\text{ compact}\subset E\}$ for all measurable $E\subset X$.} Borel probability measures and linear quasi-states on $X$. In what follows we will sometimes tacitly identify such a measure with the corresponding linear functional.
  \item It is easy to see that any quasi-state is $1$-Lipschitz with respect to the uniform norm: $|\zeta(f)-\zeta(g)|\leq\|f-g\|_{C^0}$, where $\|h\|_{C^0}=\max_{x\in X}|h(x)|$.
\end{enumerate}

\end{rem}
Quasi-states arose from a discussion around the von Neumann axioms of quantum mechanics, see \cite[Section 10]{Entov_Polterovich_QS_Sympl_intersections}, \cite[Section 1]{Butler_Semisolid_sets_TMs}, and the references therein. The above definition is due to Aarnes \cite{Aarnes_QS_QM}, who also provided a very general construction of quasi-states in terms of so-called solid sets and solid set functions \cite{Aarnes_Constr_NSA}. A subset $A\subset X$ is \tb{solid} if $A$ and $X\setminus A$ are both connected. The following result is fundamental for the present paper.
\begin{thm}\label{thm:Aarnes_constr}
  Let $X$ be a connected finite CW-complex with $H^1(X;\Z) = 0$, and let $\mu$ be a Borel probability measure on $X$ with the property that there do not exist solid closed disjoint sets $C,C'\subset X$, such that $\mu(C)=\mu(C')=\frac 12$. Then there exists a unique quasi-state $\zeta$ on $X$ with the following property: whenever $U\subset X$ is a solid open set such that $\mu(U) \leq \frac12$, then for each $f \in C(X)$ supported in $U$ we have $\zeta(f) = 0$.
\end{thm}
This result is due to a collective effort of a few people: J.\ Aarnes, F.\ Knudsen, A.\ Rustad, S.\ Butler, and others. It does not appear exactly in this form in the literature, but, at least implicitly, it has been known to the specialists for some time. See Section \ref{ss:QS_TM} for a proof and some historical remarks. For the sake of brevity and simplicity, and also because in its basic form it was proved by Aarnes in \cite[Section 6]{Aarnes_Constr_NSA}, we will refer to the above as the \tb{Aarnes construction}, and to the resulting quasi-states as \tb{Aarnes quasi-states}.

\begin{rem}
  We will only apply the Aarnes construction when the underlying space is a closed smooth manifold. This is possible since any such manifold can be given a CW structure via a triangulation \cite{Cairns_Triangulation_mfds}.
\end{rem}

The following is a basic example. It is hard to pinpoint who exactly came up with it, but it is safe to say that this was known to the specialists toward the end of the $90$'s.
\begin{exam}\label{ex:median_QS}
  In case $X=S^2$ and $\mu$ is the normalized Lebesgue measure on $X$ coming from the round area form, the Aarnes construction yields the so-called \tb{median quasi-state}. Thanks to Remark \ref{rem:Riesz_rep} (ii), $\zeta$ is uniquely determined by its values on Morse functions. If $f \in C^\infty(X)$ is Morse, the value $\zeta(f)$ can be obtained as follows. There exists a unique component $C_f$ of a level set of $f$, called its median, such that each connected component $U$ of $X\setminus C_f$ satisfies $\mu(U) \leq \frac12$. Then $\zeta(f)$ is the value of $f$ on $C_f$. A similar description can be obtained for the Aarnes quasi-state on any closed manifold $X$ and any measure on it as in Theorem \ref{thm:Aarnes_constr}.
\end{exam}

Aarnes quasi-states have a special property called \emph{simplicity}. Namely, a quasi-state $\zeta$ on $X$ is \tb{simple} if for each $f \in C(X)$ we have $\zeta(f^2)=\zeta(f)^2$, or equivalently, $\zeta$ is \emph{multiplicative} on each $C(f)$. Simple quasi-states were introduced in \cite{Aarnes_Pure_QS_Extremal_QM} (there called `pure'), and they are indecomposable, meaning they cannot be written as nontrivial convex combinations of other quasi-states, and thus they lie at the extremal boundary of the convex set of all quasi-states, see \emph{ibid.}

The following special case of the Aarnes construction will play an important role for us.
\begin{coroll}\label{cor:Aarnes_delta}
  Let $M$ be a closed connected manifold and $\mu$ a measure on $M$ as in Theorem \ref{thm:Aarnes_constr}, and assume that there is $x_0 \in M$ with $\mu(\{x_0\}) \geq \frac 12$. Then the corresponding Aarnes quasi-state is the evaluation at $x_0$.
\end{coroll}
\noindent  For completeness, we give a proof at the end of Section \ref{ss:QS_TM}. Another proof can be obtained as a consequence of \cite[Lemma 4.11]{Butler_Decomps_signed_deficient_TMs}.

The main point of the present paper is a certain triviality result for the Aarnes quasi-states in the context of symplectic geometry. We refer the reader to \cite{McDuff_Salamon_Intro} as well as Section \ref{ss:sympl_geom} for the relevant background on symplectic geometry, and to \cite{Polterovich_Rosen_Fcn_thry_sympl_mfds} for a discussion of its relationship to quasi-states. In \cite{Entov_Polterovich_QS_Sympl_intersections}, Entov--Polterovich defined a \tb{symplectic quasi-state}\footnote{In that paper the definition of a symplectic quasi-state included two additional properties, namely vanishing and symplectic invariance. However, it was later realized that these arise as a consequence of the construction, and thus they were later dropped from the definition. Here we use the more general definition which does not assume these properties.} on a closed symplectic manifold $(M,\omega)$ as a quasi-state $\zeta \fc C(M) \to \R$ satisfying the more restrictive assumption of
\begin{itemize}
  \item \tb{(strong quasi-linearity):} given Poisson-commuting $f,g \in C^\infty(M)$, $\zeta$ is linear on the subspace of $C(M)$ spanned by $f,g$.
\end{itemize}
They proved that if $\dim M = 2$, then any quasi-state on $M$ is symplectic (\emph{ibid.}, Theorem 8.1). Any linear quasi-state is trivially symplectic. The only known construction of \emph{nonlinear} symplectic quasi-states for $\dim M \geq 4$ uses Floer theory. The Floer-theoretic construction of \cite{Entov_Polterovich_QS_Sympl_intersections} in particular applies to the sphere $S^2$, where it yields the median quasi-state from Example \ref{ex:median_QS} (see \emph{ibid.}, p90, and \cite[Section 5]{Entov_Polterovich_Calabi_QM_QH}), which, being an Aarnes quasi-state, is simple. This led Entov--Polterovich to ask the following:
\begin{center}
  \emph{Are all Floer-theoretic symplectic quasi-states simple?}
\end{center}
This is \cite[Question 8.5]{Entov_Polterovich_QS_Sympl_intersections}, and although there it was only asked about $\C P^n$, the question itself, of course, makes sense for any Floer-theoretic quasi-state. Recently this has been answered in the positive as a consequence of a deep result by Mak--Sun--Varolg\"une\c s, see \cite[Corollary 1.10]{Mak_Sun_Varolgunes_Characterization_Heaviness_SH}. Since all Aarnes quasi-states are simple, as are all Floer-theoretic ones, the following natural question arises:
\begin{center}
  \emph{Can the Aarnes construction yield a symplectic quasi-state in dimension at least four?}
\end{center}
It seems that the answer should be ``obviously no,'' since the measure used in the Aarnes construction has no relation to the symplectic form, however, as we will see, the situation is far from trivial, see Remark \ref{rem:Liouville_measure_still_not_sympl}. Note that Corollary \ref{cor:Aarnes_delta} implies that the Aarnes construction does yield a symplectic quasi-state when there is a point of measure $\geq\frac 12$, but it is for a trivial reason, because the resulting quasi-state is linear.

Our main result then reads as follows.
\begin{thm}[\tb{Main result}]\label{thm:main_result}
  Let $(M,\omega)$ be a closed connected symplectic manifold of dimension at least four, with $H^1(M;\Z)=0$, and let $\zeta$ be an Aarnes quasi-state on $M$. If $\zeta$ is symplectic, then it is a delta-measure.
\end{thm}

\begin{rem}\label{rem:Liouville_measure_still_not_sympl}
  Note that even when the measure in question is the normalized Lebesgue measure corresponding to the volume form $\omega^{\wedge\frac12\dim M}$, that is the measure `remembers' some information encoded in the symplectic form, the corresponding Aarnes quasi-state is still not symplectic. This conforms to the deep distinction between symplectic geometry and the geometry of a volume form, which was first established in Gromov's seminal paper \cite{Gromov_Pseudoholomorphic}.
\end{rem}

\subsection{Symplectic embeddings and an overview of the proof}\label{ss:sympl_emb_overview_pf}

It is at this point that symplectic embeddings enter our story, in the following manner. Since we are assuming $\zeta$ to be symplectic, it is natural to try and investigate its interaction with Poisson-commuting functions on $M$, which necessitates the construction of such functions. It is easy to do this locally on any symplectic manifold, thanks to Darboux's theorem: indeed, any point in a symplectic manifold has a neighborhood which is symplectomorphic to a Euclidean ball $B\subset \C^n$ with the standard symplectic form. The functions $|z_1|^2,\dots,|z_n|^2$ Poisson commute, and we can use appropriate cutoffs to produce Poisson commuting functions supported in $B$, which can then be extended by zero to $M$, while preserving the commutation.

The issue is that $B$ may be small relative to the measure used in the construction of $\zeta$,  and thus $\zeta$ will vanish on functions supported in $B$. We are thus led to the problem of embedding `standard' sets into a given symplectic manifold so as to cover a large portion of the given measure. The following describes a family of standard sets which is crucial for us.
\begin{notation}\label{not:disks_polydisks}
  For $r>0$ we let $D(r)=\{z \in \C\,|\,|z|^2<r\}$ be the corresponding open complex disk. A \tb{symplectic polydisk} (or simply a \tb{polydisk}) is a subset of $\C^n$ of the form $P(a)=\prod_{j=1}^n D(a_j)$, where $a=(a_1,\dots,a_n)\in(0,\infty)^n$.
\end{notation}
The relevant symplectic embedding result is as follows.
\begin{thm}\label{thm:main_embedding}
  Let $(M,\omega)$ be a closed connected symplectic manifold, and let $\mu$ be a Borel probability measure on $M$. Then for each $\epsilon>0$ there exists a symplectic polydisk $P$ and a symplectic embedding $\iota \fc P \hookrightarrow M$ such that $\mu(M\setminus\iota(P)) < \epsilon$.
\end{thm}

\begin{rem}
\begin{itemize}
  \item In case $\mu$ is the normalized Lebesgue measure corresponding to the volume form $\omega^{\wedge\frac12\dim M}$, this was proved by Schlenk for any connected symplectic manifold of finite volume, see \cite[Theorem 6.1.1]{Schlenk_Embeddings_problems}. We believe our proof can be adapted to the case when $M$ is open, that is non-compact and with empty boundary, however as stated, the embedding result suffices to prove our main result.
  \item Just like in \cite{Schlenk_Embeddings_problems}, we use the technique of folding to construct our embedding, and use polydisks of the form $P(\alpha,\eta,\dots,\eta)$. The issue of optimality of the ratio $\frac{\alpha}{\eta}$ of the parameters of these polydisks, given the measure $\mu$, seemingly belongs to the realm of `hard' symplectic topology, and is not treated here.
\end{itemize}

\end{rem}

To see how this embedding result helps us prove Theorem \ref{thm:main_result}, we need to discuss involutive maps.
\begin{defin}\label{def:involutive_map}
  Let $(M,\omega)$ be a symplectic manifold. A smooth map $\Phi \fc M \to \R^k$ is \tb{involutive} if its coordinate functions pairwise Poisson commute.
\end{defin}

The significance of involutive maps lies in their relation with symplectic quasi-states. To formulate it, let us recall the notion of pushforward.\footnote{In the present context this appears, for instance in \cite{Rustad_QMs_with_image_transformations}, where it is called the adjoint of an image transformation. See also reference therein.} If $X,Y$ are compact Hausdorff spaces, $\Phi\fc X \to Y$ is a continuous map, and $\zeta$ is a quasi-state on $X$, then the functional $\Phi_*\zeta \fc C(Y) \to\R$ defined for $g \in C(Y)$ by $\Phi_*\zeta(g)=\zeta(g\circ \Phi)$ is again a quasi-state, called the \tb{pushforward of $\zeta$ by $\Phi$}. The following lemma is proved at the end of Section \ref{ss:sympl_geom}.
\begin{lemma}\label{lem:char_sympl_QS_invol_maps}
  Let $(M,\omega)$ be a closed symplectic manifold, and let $\zeta$ be a quasi-state on $M$. Then $\zeta$ is symplectic if and only if for each involutive map $\Phi \fc M \to \R^k$ the functional $\Phi_*\zeta$, viewed as a quasi-state on the compact space $\im\Phi$, is linear.
\end{lemma}

Now assume that $(M,\omega)$ and $\mu$ are as in Theorem \ref{thm:main_result}, and that the resulting Aarnes quasi-state $\zeta$ is symplectic. Using an embedding $\iota \fc P \hookrightarrow M$ as in Theorem \ref{thm:main_embedding}, we can produce an involutive map $\Phi \fc M \to \R^n$, where $n=\frac12\dim M$, by taking the functions $|z_j|^2$ on $P$, suitably cutting them off so that they are supported in $P$, and implanting the result into $M$ via $\iota$.

The idea now is that since $\zeta$ is symplectic and $\Phi$ is involutive, $\Phi_*\zeta$ is linear by Lemma \ref{lem:char_sympl_QS_invol_maps}, and thus, by Remark \ref{rem:Riesz_rep}, it is the integral with respect to a regular Borel probability measure on $\im\Phi$. Since $\zeta$ is simple, so is $\Phi_*\zeta$, and since the only regular Borel probability measures representing simple linear quasi-states are delta-measures (see Lemma \ref{lem:simple_meas_delta}), we see that $\Phi_*\zeta$ is the evaluation at some point $x_0\in\im\Phi$.

Choosing $\Phi$ carefully, we can use this to construct a Lagrangian torus $T\subset M$, containing a suitable connected component of $\Phi^{-1}(x_0)$, with the property that if $f\in C(M)$ satisfies $f|_T\equiv c$, then $\zeta(f) = c$. Following \cite{Entov_Polterovich_Rigid_subsets}, we call sets with this property \tb{$\zeta$-superheavy}. The last piece of the puzzle is the following theorem, which, together with the embedding result, Theorem \ref{thm:main_embedding}, forms the technical heart of the present paper.
\begin{thm}\label{thm:sh_Lag_torus_contains_sh_point}
  Let $(M,\omega)$ be a closed connected symplectic manifold of dimension at least four, with $H^1(M;\Z)=0$, and let $\zeta$ be the Aarnes quasi-state on $M$ corresponding to a measure $\mu$ as in Theorem \ref{thm:Aarnes_constr}, and assume that $\zeta$ is symplectic. If $T\subset M$ is a $\zeta$-superheavy Lagrangian torus, then there is $z_0 \in T$ with $\mu(\{z_0\}) \geq \frac 12$.
\end{thm}
\noindent Corollary \ref{cor:Aarnes_delta} then implies that $\zeta$ is the evaluation at $z_0$. This result, in combination with the above discussion, implies our main result, Theorem \ref{thm:main_result}. Theorem \ref{thm:sh_Lag_torus_contains_sh_point} is proved in Section \ref{ss:sh_Lag_has_sh_point}.

\begin{rem}
  The construction of a superheavy Lagrangian torus $T\subset M$ under the assumptions of Theorem \ref{thm:main_result} and the construction of a point $z_0 \in T$ with $\mu(\{z_0\}) \geq \frac12$ are logically independent, therefore we separated these into two distinct statements.
\end{rem}

\subsection{Discussion}\label{ss:discussion}

The topic of constraints on symplectic quasi-states started with the following unpublished result due to L.\ Polterovich, which he proved circa 2007:
\begin{thm}\label{thm:Polterovich_restriction}
  Let $M=S^2\times S^2$ be endowed with the symplectic form $\omega\oplus\omega$, where $\omega$ is the standard round area form. Let $\zeta$ be a symplectic quasi-state on $M$, and denote by $\pi_{1,2}\fc M \to S^2$ the projections to the factors. Then there do not exist triples of pairwise distinct points $a_i,b_i,c_i \in S^2$, $i=1,2$, that is $a_i\neq b_i\neq c_i\neq a_i$, such that for $i=1,2$, $\pi_{i*}\zeta$ is the Aarnes quasi-state corresponding to the measure $\frac13(\delta_{a_i}+\delta_{b_i}+\delta_{c_i})$.
\end{thm}
\noindent The basic idea of his proof is the use of the involutive maps coming from two distinct integrable systems. Compare this to the proof of our main result, which uses infinitely many different involutive maps, see Section \ref{s:pf_main_result}.

\begin{rem}
  If in the formulation of Theorem \ref{thm:Polterovich_restriction} we replace the Aarnes quasi-states corresponding to measures of the form $\frac13(\delta_a+\delta_b+\delta_c)$ by the median quasi-state from Example \ref{ex:median_QS}, then the conclusion is no longer true, since it can be easily shown that if $\zeta$ is the Entov--Polterovich symplectic quasi-state on $S^2\times S^2$, then its projections $\pi_{i*}\zeta$ both equal the median quasi-state.
\end{rem}

Certain constraints on Aarnes quasi-states were also proved in \cite{Dickstein_Zapolsky_Approx_QS_manifolds}. \emph{Ibid.}, we introduced, following a suggestion of L.\ Polterovich, the Wasserstein distance\footnote{In the context of transportation theory, this is also known as the Kantorovich--Rubinstein distance.} $W_1$ on the space of quasi-states on a given compact metric space $(X,d)$, namely
$$W_1(\zeta,\eta) = \max\{\zeta(f)-\eta(f)\,|\,f\in \Lip_1(X,d)\}\,,$$
where $\Lip_1(X,d)$ is the collection of Lipschitz functions on $X$ with Lipschitz constant at most $1$. In the same paper, Corollary 1.14 states that the $W_1$-distance between any Entov--Polterovich symplectic quasi-state on $M=\C P^n$ or $M=S^2\times S^2$ and the collection of Aarnes quasi-states on $M$ is positive.

Returning to Theorem \ref{thm:Polterovich_restriction}, if we allow two of the three points $a_i,b_i,c_i$ to coincide, say, $x_i=a_i=b_i$ for $i=1,2$, then the corresponding Aarnes quasi-states are the delta-measures at $x_i$, thanks to Corollary \ref{cor:Aarnes_delta}, and they, of course, \emph{can} be obtained by projecting a symplectic quasi-state on $S^2\times S^2$, namely the delta-measure at $(x_1,x_2)$. This hints at the existence of a quantitative constraint for a symplectic quasi-state $\zeta$ on $S^2\times S^2$, namely that the $W_1$-distance from its projections $\pi_{i*}\zeta$ to such three-point Aarnes quasi-states should be governed by the $W_1$-distance from the Aarnes quasi-states to the collection of delta-measures. This is indeed the case, as is illustrated by the following result, which will appear in \cite{Dickstein_PhD_2025}. To formulate it, let $\cS(M,\omega)$ and $\Delta(X)$ stand for the collections of symplectic quasi-states on $(M,\omega)$ and of delta-measures on $X$, respectively. Also, let $\cT(S^2)$ be the collection of Aarnes quasi-states on $S^2$ obtained from measures of the form $\frac13(\delta_a+\delta_b+\delta_c)$ for three pairwise distinct points $a,b,c \in S^2$. Then we have
\begin{thm}
  $$\inf_{\substack{{\eta_1,\eta_2 \in \cT(S^2)} \\ {\zeta \in \cS(S^2\times S^2,\omega\oplus\omega)}}} \left(\frac{W_1(\eta_1,\pi_{1*}\zeta)}{W_1(\eta_1,\Delta(S^2))} + \frac{W_1(\eta_2,\pi_{2*}\zeta)}{W_1(\eta_2,\Delta(S^2))}\right) > 0\,.$$
\end{thm}

Let us relate this discussion to our main result, Theorem \ref{thm:main_result}. It says that the intersection of the collections of Aarnes quasi-states and of symplectic quasi-states on a closed connected symplectic manifold $M$ with $H^1(M;\Z)=0$ and $\dim M\geq 4$ consists exactly of delta-measures. It would be interesting to prove a quantitative result of this sort. We propose the following
\begin{conj}
  Endow $M$ with a metric inducing its topology. Then there exists a neighborhood $\cU$ of the set of delta-measures in the space of quasi-states on $M$, and a constant $c \in (0,1]$ such that for any Aarnes quasi-state $\zeta \in \cU$ we have
  $$W_1(\zeta,\cS(M,\omega)) \geq c\cdot W_1(\zeta,\Delta(M))\,.$$
\end{conj}

Lastly, we would like to point out that there are more general `soft' constructions of quasi-states. In \cite{Aarnes_Constr_NSA} Aarnes defined a certain quantity associated to a connected locally connected compact Hausdorff topological space, called \emph{genus}, see \emph{ibid.}, the paragraph immediately after Definition 2.2 on page 215, where it is denoted by $g$. In \cite{Knudsen_Topology_constr_extreme_TM} Knudsen showed that a finite CW complex $X$ has Aarnes genus zero if and only if $H^1(X;\Z)=0$. On spaces of Aarnes genus zero, there is a construction in \cite{Butler_Q_fcns} in terms of q-functions. On spaces of Aarnes genus $1$, which, for instance, include all the tori $\T^n$, Grubb describes in \cite{Grubb_Irr_partitions_constr_QMs} a general construction of quasi-states. Finally, he recently communicated to us a construction of nonlinear quasi-states, including simple ones, on very general spaces, which in particular includes all closed connected manifolds with $H^1(M;\Z)\neq 0$. An earlier construction for closed connected orientable surfaces of all genera appears in \cite{Zapolsky_Isotopy_invt_TMs_surfaces}.

A natural question therefore is whether any of these constructions yield nonlinear symplectic quasi-states in case the symplectic manifold in question has dimension $\geq 4$. In view of general principles of symplectic topology, as well as our main result Theorem \ref{thm:main_result}, a positive answer seems to us to be extremely unlikely, however a proof of this in any of the aforementioned cases will undoubtedly be very technically challenging. This is already illustrated by the complexity of the arguments in Section \ref{ss:sh_Lag_has_sh_point} involved in the proof of our main result in the relatively simple setting of Aarnes quasi-states.

\paragraph*{Organization of the paper.} We collect the necessary preliminaries on measure theory, symplectic geometry, and quasi-states in Section \ref{s:prelims}, in particular in Section \ref{ss:QS_TM} we present the proofs of Theorem \ref{thm:Aarnes_constr} and Corollary \ref{cor:Aarnes_delta}, as well as some historical remarks. The proof of the embedding result, Theorem \ref{thm:main_embedding}, occupies Section \ref{s:constructing_emb}. Section \ref{s:pf_main_result} is dedicated to the proof of our main result, Theorem \ref{thm:main_result}, which is split into two parts: the first one, described in Section \ref{ss:involutive_map_fibers}, uses our embedding result to construct an involutive map $\Phi \fc M \to \R^n$, one of whose fibers has a component contained in a superheavy Lagrangian torus. The second part, Section \ref{ss:sh_Lag_has_sh_point}, is dedicated to the proof of Theorem \ref{thm:sh_Lag_torus_contains_sh_point}, which happens entirely in a small Weinstein neighborhood of the superheavy Lagrangian torus.

\paragraph*{Acknowledgements.} We wish to thank Leonid Polterovich for valuable discussions and in particular for pointing out Schlenk's book \cite{Schlenk_Embeddings_problems} to us. We are grateful to Svetlana Butler for an illuminating discussion about the history of the field of quasi-states, as well as references and useful remarks. We extend our gratitude to the anonymous referee for carefully reading the paper, making poignant remarks regarding some points in the text, and especially for providing us with a very elegant and much simpler proof of Theorem \ref{thm:sh_Lag_torus_contains_sh_point} than we originally had, shortening the text by a good few pages. Most of the paper was written when the first author was a doctoral student at Tel Aviv University, while some of the work was done when he was a postdoctoral fellow at Universit\'e de Montr\'eal, and he wishes to acknowledge its hospitality and excellent research atmosphere.

\section{Preliminaries}\label{s:prelims}

Here we gather the necessary fundamental notions and results, which will be used in the rest of the paper. Throughout we will use the $\ell^\infty$-norm on $\R^N$: $\|x\|=\max_{1\leq i\leq N}|x_i|$, and its associated metric $d(x,y) = \|x-y\|$ for $x,y \in \R^N$. We let $B_x(r)$, $\ol B_x(r)$ be the corresponding open, respectively, closed balls of radius $r$ centered at $x \in \R^N$. We note that
$$B_x(r) = \prod_{i=1}^N (x_i-r,x_i+r)\quad \text{and}\quad \ol B_x(r) = \prod_{i=1}^N [x_i-r,x_i+r]\,.$$
That is, the $\ell^\infty$-balls are in fact cubes. We also define the translation operator corresponding to $y \in \R^N$ by $T_y \fc \R^N\to \R^N$, $T_y(x) = x+y$.

\subsection{Measure theory}\label{ss:measure_thry}

The following is well-known in measure theory, but we give a proof for the sake of completeness.
\begin{lemma}\label{lem:pairwise_disjoint_meas_sets_at_most_countable}
  In a measure space of finite total measure, any collection of pairwise disjoint measurable sets, all of which have positive measure, is at most countably infinite.
\end{lemma}
\begin{proof}
  Let $\cA$ be such a collection. For $n \in \N$ put $\cA_n = \{A \in \cA\,|\,A\text{ has measure }>\frac 1 n\}$. Then $\cA = \bigcup_{n \in \N}\cA_n$. By assumption, each $\cA_n$ is a finite collection and the assertion follows.
\end{proof}
\noindent We will also use the following lemmas.
\begin{lemma}[{\cite[Theorem 1.19 (d-e)]{Rudin_Real_complex_analysis}}]\label{lem:decreasing_seq_meas_sets_limit}
  In a measure space of finite total measure, for a decreasing sequence $(A_i)_{i\in\N}$ of measurable sets we have $\mu\big(\bigcap_{i\in\N}A_i\big) = \lim_{i\to\infty}\mu(A_i)$. If $(A_i)_{i\in\N}$ is an increasing sequence of measurable sets, then $\mu\big(\bigcup_{i\in\N}A_i) = \lim_{i\to\infty}\mu(A_i)$. \qed
\end{lemma}

\begin{lemma}\label{cor:approx_open_set_cpts}
  Let $U\subset \R^N$ be an open bounded set, and let $\nu$ be a finite Borel measure on $\R^N$. Then
  $$\nu(U) = \sup\{\nu(K)\,|\,K\text{ compact}\subset U\}\,. \qquad \qed$$
\end{lemma}
\noindent This follows from \cite[Chapter 6, Proposition 1.3]{Stein_Shakarchi_Real_analysis_3}.

\begin{figure}
  \centering
  \includegraphics[width=8cm]{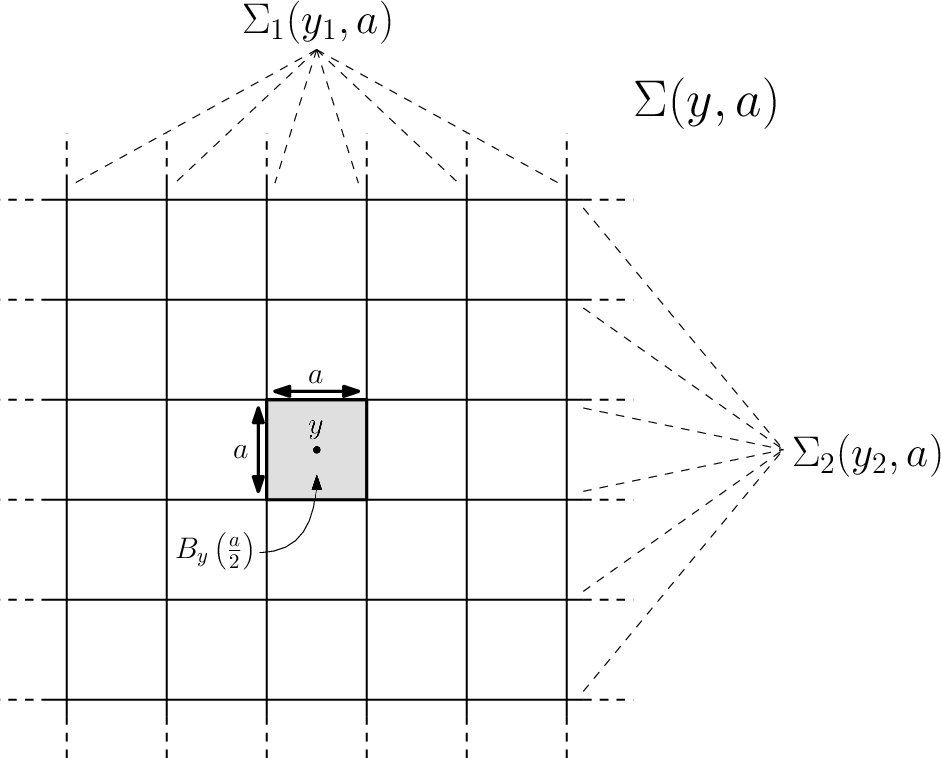}
  \caption{The lattice $\Sigma(y,a)$ in dimension $2$}\label{fig:Sigma_lattice}
\end{figure}

For $i=1,\dots,N$, $u \in \R$, $y \in \R^N$, and $a>0$ we let
$$\Sigma_i(u,a) = \{x \in \R^N\,|\,x_i\equiv \tfrac a 2 + u\bmod a\Z\}\quad\text{and}\quad \Sigma(y,a) = \bigcup_{i=1}^N\Sigma_i(y_i,a)\,.$$
See Figure \ref{fig:Sigma_lattice}. For future use let us also define
$$\textstyle\cU(y,a) = \big\{B_x(\tfrac a 2)\,|\,\forall i: x_i \equiv y_i \bmod a\Z\big\}\,,\quad\cC(y,a) = \big\{\ol B{}_x(\tfrac a 2)\,|\,\forall i: x_i \equiv y_i \bmod a\Z\big\}\,,$$
that is $\cU(y,a)$ is the set of connected components of $\R^N\setminus\Sigma(y,a)$ while $\cC(y,a)$ is the set of their closures. We will use the following basic property.

\begin{lemma}\label{lem:Sigma_y_a_at_most_countable}
  Let $\nu$ be a finite Borel measure on $\R^N$. Then for any $a > 0$ the collection of $y \in \R^N$ such that $\nu(\Sigma(y,a))\neq 0$ has Lebesgue measure zero. In particular, for each open $U\subset\R^N$ there is $y \in U$ such that $\nu(\Sigma(y,a))=0$.
\end{lemma}
\begin{proof}
  Fix $a>0$ and $i=1,\dots,N$ and consider the collection $\{\Sigma_i(b,a)\}_{b\in[0,a)}$. It consists of pairwise disjoint measurable subsets, therefore by Lemma \ref{lem:pairwise_disjoint_meas_sets_at_most_countable}, the set
  $$\cB_i:=\{b \in [0,a)\,|\,\nu(\Sigma_i(b,a)) \neq 0\}$$
  is at most countable. We have $\{b \in \R\,|\,\nu(\Sigma_i(b,a))\neq 0\}=\cB_i+a\Z$, which is likewise at most countable. Since a finite union of measurable sets has measure zero if and only if so does each set in the union, and since $\Sigma(y,a) = \bigcup_{i=1}^N \Sigma_i(y_i,a)$, we see that $\nu(\Sigma(y,a)) = 0$ if and only if for each $i=1,\dots,N$ we have $\nu(\Sigma_i(y_i,a))=0$, whence
  $$\{y\in\R^N\,|\,\nu(\Sigma(y,a)) = 0\} = \bigcap_{i=1}^N\{y\in\R^N\,|\,\nu(\Sigma_i(y_i,a))=0\}=\prod_{i=1}^{N}\R\setminus(\cB_i+a\Z)\,.$$
  Letting $\pr_i\fc\R^N\to\R$ be the projection to the $i$-th coordinate, we see that
  \begin{multline*}
    \{y\in\R^N\,|\,\nu(\Sigma(y,a)) \neq 0\} = \R^N\setminus\{y\in\R^N\,|\,\nu(\Sigma(y,a)) = 0\}\\
      =\R^N\setminus\prod_{i=1}^{N}\R\setminus(\cB_i+a\Z) = \bigcup_{i=1}^N\pr_i^{-1}(\cB_i+a\Z)\,.
  \end{multline*}
  Since for each $i$, $\cB_i+a\Z$ is at most countable, the preimage $\pr_i^{-1}(\cB_i+a\Z)$ has Lebesgue measure zero in $\R^N$, and therefore so does $\{y\in\R^N\,|\,\nu(\Sigma(y,a)) \neq 0\}$ as a finite union of such sets.
\end{proof}

\subsection{Symplectic geometry}\label{ss:sympl_geom}

The material here is standard, but we include it for the sake of readers who are not symplectic geometers, and to fix notation. Some of the examples presented here will be necessary in the next sections.

A \tb{symplectic manifold} is a pair $(M,\omega)$, where $M$ is an even-dimensional smooth manifold, while $\omega$ is a symplectic form on $M$, that is $\omega$ is a closed $2$-form, which is nondegenerate, that is $\omega^{\wedge \frac12\dim M}$ is a volume form. Smooth functions on $M$ are also referred to as \tb{Hamiltonians}, and given $f \in C^\infty(M)$ its \tb{Hamiltonian vector field} $X_f$ is determined by $\iota_{X_f}\omega=-df$. The \tb{Poisson bracket} of $f,g \in C^\infty(M)$ is the function $\{f,g\} = \omega(X_g,X_f) = df(X_g)$.

If the Hamiltonian vector field $X_f$ of $f$ is complete, it integrates into the corresponding \tb{Hamiltonian isotopy} $\phi_f^t$ such that $\phi_f^0=\id_M$. We denote $\phi_f=\phi_f^1$. A \tb{symplectomorphism} of $(M,\omega)$ is a diffeomorphism preserving $\omega$. A crucial feature of symplectic geometry is that Hamiltonian isotopies preserve $\omega$, and thus smooth functions yield a plethora of symplectomorphisms, also known as \tb{Hamiltonian diffeomorphisms}. More generally, if $(M,\omega)$ and $(M',\omega')$ are symplectic manifolds, a diffeomorphism $\phi \fc M \to M'$ is called a symplectomorphism if $\phi^*\omega'=\omega$. If such $\phi$ exists, we say that $(M,\omega)$, $(M',\omega')$ are \tb{symplectomorphic}.

\begin{exam}\label{ex:sympl_mfds}
  \begin{itemize}
    \item The Euclidean space $\R^{2n}(x_1,y_1,\dots,x_n,y_n)$ carries the standard symplectic form $\omega_{\std}=\sum_{j=1}^{n}dx_j\wedge dy_j$. We tacitly identify $\R^{2n}=\C^n$ via $(x,y)\leftrightarrow z=x+iy$.
    \item If $(M,\omega)$ is a symplectic manifold, then any open subset $U\subset M$ becomes a symplectic manifold once equipped with the restricted form $\omega|_U$. \emph{Darboux's theorem}, already mentioned above, states that any $2n$-dimensional symplectic manifold $(M,\omega)$ is locally symplectomorphic to $(\R^{2n},\omega_{\std})$, meaning any point of $M$ has an open neighborhood $U$ such that $(U,\omega|_U)$ is symplectomorphic to an open subset of $\R^{2n}$.
    \item If $Q$ is a smooth manifold of dimension $n$, its cotangent bundle $T^*Q$ carries the canonical symplectic form $\omega^Q=d\lambda^Q$, where $\lambda^Q$ is the Liouville $1$-form, given in local coordinates $(q,p)$ on $T^*Q$ by $\lambda^Q=\sum_{j=1}^{n}p_i\,dq_i$.
    \item A surface with an area form is symplectic, since the area form is automatically closed.
    \item The product of two symplectic manifolds $(M,\omega)$, $(M',\omega')$ is again symplectic with the form $\omega\oplus\omega'$. For instance, $\R^{2n}$ is the product of $n$ copies of $(\R^2(x,y),dx\wedge dy)$.
  \end{itemize}
\end{exam}

Next we will describe some basic examples of symplectomorphisms between certain open subsets of $\R^{2n}$. The reader is also invited to consult \cite{Schlenk_Embeddings_problems}, especially Section 1.3 there.

\begin{exam}\label{ex:symplectomorphisms}
  \begin{itemize}
    \item \tb{Translations.} Any translation operator on $\R^{2n}$ is a symplectomorphism.
    \item \tb{Area-preserving diffeomorphisms.} Since in dimension two a symplectic form is simply an area form, two open subsets of $\R^2$ are symplectomorphic if and only if they are diffeomorphic and have the same area, see, for instance \cite[Proposition 1]{Schlenk_Embeddings_problems}. Thus, for $a>0$, the open disk $D(a)$ is symplectomorphic to the square $(0,\sqrt{\pi a})^2$.
    \item \tb{Products.} The product of two symplectomorphisms is again a symplectomorphism. A typical example is a symplectomorphism between a polydisk $P(a) = \prod_{j=1}^{n}D(a_j)$ and the box $\prod_{j=1}^{n}(0,\sqrt{\pi a_j})^2$.
    \item \tb{Symplectic shears and smears.} If $f \fc \R \to \R$ is a smooth function, the corresponding \tb{symplectic shear} $S_f \fc \R^2 \to \R^2$ is defined by $S_f(x,y) = (x,y+f(x))$, and it is a symplectomorphism. If $g \fc \R \to (0,\infty)$ is a smooth function, then the corresponding \tb{symplectic smear} is the map
        $$\Theta_g \fc \R^2\to \R^2\quad \text{given by} \quad \Theta_g(x,y) = \big(\textstyle\int_0^x\frac{dt}{g(t)},g(x)y\big)\,.$$
        It is a symplectomorphism.
    \item \tb{Restrictions.} If $\phi \fc M \to M'$ is a symplectomorphism and $U\subset M$ is an open subset, then $\phi|_U\fc U \to \phi(U)$ is again a symplectomorphism.
    \item \tb{Symplectic lifts.} Given a smooth function $\rho \fc \R \to \R$, we can define a class of symplectomorphisms of $\R^{2n}$, known as \tb{symplectic lifts}. These are Hamiltonian diffeomorphisms coming from `lifting Hamiltonians' $H\in C^\infty(\R^{2n})$, defined by an expression of the form
$$H(x_1,y_1,\dots,x_n,y_n) = \rho(x_1)(ax_j+b)$$
for some numbers $a,b\in\R$, where $j\geq 2$. We can compute:
$$\phi_H(x,y) = (x,y) + \rho'(x_1)(ax_j+b)\partial_{y_1} + a\rho(x_1)\partial_{y_j}\,.$$
Similarly we can use lifting Hamiltonians of the form
$$G(x_1,y_1,\dots,x_n,y_n) = \rho(x_1)(ay_j+b)\,,$$
for which we have
$$\phi_G(x,y) = (x,y) + \rho'(x_1)(ay_j+b)\partial_{y_1} - a\rho(x_1)\partial_{x_j}\,.$$
A typical situation is when $\rho(t)\equiv 0$ for $t\leq t_0$ and $\rho(t)\equiv 1$ for $t\geq t_1$ for some $t_0<t_1$. Then $\phi_H,\phi_G$ are the identity for $x_1\leq t_0$, while for $x_1\geq t_1$, $\phi_H$ is the shift $(x,y)\mapsto (x,y) + a\partial_{y_j}$, and $\phi_G$ is the shift $(x,y) \mapsto (x,y) - a\partial_{x_j}$.
  \end{itemize}
\end{exam}

\begin{rem}\label{rem:smears}
  We will use certain elementary properties of symplectic smears. Namely, assume $h \fc \R \to \R$ is a smooth positive function. Then $x\mapsto \int_0^xh(t)\,dt$ is a diffeomorphism of $\R$; let $f$ be the its inverse, that is $f^{-1}(x) = \int_0^x h(t)\,dt$ for all $x \in \R$, and let $g=1/f'$. Then an easy calculation shows that for each $r \in \R$, the smear $\Theta_g$ maps the line $\R \times \{r\}$ to the graph of $rh$. Moreover, if $a<b$ are real numbers such that $h|_{[a,b]}\equiv 1$, then the restriction of $\Theta_g$ to the strip $\big[\int_0^ah(t)\,dt\,, b-a+\int_0^ah(t)\,dt\big] \times \R$ coincides with the shift operator $T_{(a-\int_0^ah(t)\,dt,0)}$. Another easy observation is that smears map vertical lines linearly into vertical lines. It follows that for $a > 0$, the image under $\Theta_g$ of a strip $\big(0,\int_0^ah(t)\,dt\big)\times(0,1)$ is the open band bounded by $x=0$ on the left, $x=a$ on the right, $y=0$ on the bottom and the graph of $h$ on the top.
\end{rem}

We now define one of the fundamental concepts for this paper. Given symplectic manifolds $(U,\eta)$, $(M,\omega)$ of the same dimension, a smooth map $\iota \fc U\to M$ is called a \tb{symplectic embedding}, if it is a smooth embedding, that is a diffeomorphism onto the image, such that in addition $\iota^*\omega=\eta$. Constructing symplectic embeddings and obstructing their existence is one of the central themes in symplectic topology, see \cite{Schlenk_Sympl_emb_Old_new}. An obvious example of a symplectic embedding is the inclusion of an open subset of a symplectic manifold. Also note that symplectomorphisms are precisely the invertible symplectic embeddings.

\begin{rem}\label{rem:injective_sympl}
  Keeping the notation of the previous paragraph, let us call a smooth map $\iota \fc U \to M$ symplectic if $\iota^*\omega=\eta$. In this case $\iota^*(\omega^{\wedge n})=\eta^{\wedge n}$, where $n=\frac12\dim M$, and thus $\iota$ is an immersion, therefore it is a local diffeomorphism. If in addition it is injective, then it is an embedding. Thus symplectic embeddings are simply injective symplectic maps.
\end{rem}

Our main embedding result, Theorem \ref{thm:main_embedding}, is based on an interplay between two shapes---symplectic polydisks, see Notation \ref{not:disks_polydisks}---and \emph{cubes}:
\begin{notation}\label{not:cube}
  The \tb{symplectic cube} $C^{2n}(a)\subset \R^{2n}$ with side $a > 0$, is the Euclidean cube $(0,a)^{2n}$ endowed with the restriction of $\omega_{\std}$.
\end{notation}

Above we have defined involutive maps, see Definition \ref{def:involutive_map}. Let us give some examples.
\begin{exam}\label{ex:invol_maps}
  \begin{enumerate}
    \item If $(M,\omega)$ is a symplectic manifold, then any function $f \fc M \to \R$ is trivially an involutive map.
    \item If $(M,\omega)$ and $(M',\omega')$ are symplectic manifolds and $\Phi \fc M \to \R^n$, $\Phi' \fc M' \to \R^{n'}$ are involutive maps, then so is the product $\Phi \times \Phi' \fc M \times M' \to \R^{n+n'} = \R^n\times\R^{n'}$.
    \item A special case is as follows: if $(M_i,\omega_i)$ are symplectic manifolds for $i=1,\dots,k$, and $f_i\in C^\infty(M_i)$, then $f_1\times\dots\times f_k \fc M_1\times\dots\times M_k\to\R^k$ is involutive.
    \item Combining the previous item with the observation that $\C^n$ is, as a symplectic manifold, the product of $n$ copies of $(\C,\omega_{\std})$, we obtain the involutive map $\Phi_{\std} \fc \C^n \to \R^n$, $\Phi_{\std}(z) = (|z_1|^2,\dots,|z_n|^2)$. Of course, we can use any function of $z_j$ in the $j$-th coordinate, but this specific involutive map will play a crucial role in the proof of the main result, which is why we have singled it out.
    \item Another example is given as follows. If we identify $S^1=\R/\Z$, then the cotangent bundle $T^*S^1$ can be identified with $S^1(q)\times\R(p) = \R/\Z\times\R$ with the symplectic form given by $\omega^{S^1}=dp\wedge dq$. Of course, the projection to the second factor $T^*S^1=\R/\Z\times\R\to\R$, viewed as a function, is an involutive map. Taking the product of $n$ copies of this, we arrive at the involutive map $\pr \fc (T^*S^1)^n\to\R^n$. Relative to the identifications $(T^*S^1)^n=T^*((S^1)^n)=T^*\T^n=\T^n\times\R^n$, $\pr$ is simply the projection to the second factor. This involutive map will play a crucial role in the proof of Theorem \ref{thm:sh_Lag_torus_contains_sh_point}.
    \item If $\Phi \fc M \to \R^k$ is involutive, then so is $\Phi \circ \psi$ for any symplectomorphism $\psi \fc N \to M$.
    \item Finally, if $\Phi \fc M \to \R^k$ is an involutive map, and $\phi \fc \R^k \to \R^\ell$ is any smooth map, then $\phi \circ\Phi\fc M \to \R^\ell$ is also involutive.
  \end{enumerate}
\end{exam}

Let us discuss one final notion from symplectic geometry. A submanifold $L$ of a symplectic manifold $(M,\omega)$ is \tb{isotropic} if $\omega|_L\equiv 0$. Due to the nondegeneracy of $\omega$, this forces $\dim L \leq \frac12\dim M$. The boundary case $\dim L=\frac12\dim M$ is fundamental to the whole of symplectic geometry; in this case $L$ is called \tb{Lagrangian}.
\begin{rem}\label{rem:Weinstein_nbd_thm}
For Lagrangian submanifolds we have the \emph{Weinstein neighborhood theorem} \cite[Theorem 3.4.13]{McDuff_Salamon_Intro}: If $L$ is a closed Lagrangian submanifold, then there is an open neighborhood of $L$ which is symplectomorphic to an open neighborhood of the zero section in $T^*L$, endowed with the canonical form $\omega^L$, such that $L$ maps to the zero section by the identity map. In the particular case when $L$ is a Lagrangian torus, we conclude that there is an open neighborhood of $L$ symplectomorphic to an open neighborhood of $\T^n\times\{0\}$ in $T^*\T^n=\T^n\times\R^n$, where $n=\frac12\dim M$.
\end{rem}

\begin{exam}\label{ex:Lag_tori_T_alpha}
  Consider the involutive map $\Phi_{\std}\fc\C^n\to\R^n$ from Example \ref{ex:invol_maps}, item (iv). Its image is the positive orthant $[0,\infty)^n$. For $\alpha \in [0,\infty)^n$ the set
  $$T(\alpha):=\Phi_{\std}^{-1}(\alpha) = \{z \in \C^n\,|\,\forall j:|z_j|^2=\alpha_j\}$$
  is an isotropic torus of dimension equal to the number of nonzero $\alpha_j$. Thus for $\alpha \in (0,\infty)^n$, $T(\alpha)$ is Lagrangian.
\end{exam}

We close this subsection with
\begin{proof}[Proof of Lemma \ref{lem:char_sympl_QS_invol_maps}]Assume that $\zeta$ is symplectic, and let $\Phi \fc M \to \R^k$ be involutive. If $f,g \in C(\im\Phi)$, we will show that $\Phi_*\zeta$ is linear on $\spn_\R(f,g)$. Consider sequences $f_k,g_k\in C^\infty(\R^k)$ such that $f_k|_{\im\Phi}\to f$, $g_k|_{\im\Phi} \to g$ uniformly. Since $\Phi$ is involutive, for each $k$ the functions $\Phi^*f_k$, $\Phi^*g_k$ Poisson commute, and thus for $a,b\in\R$ we have $\zeta(a\Phi^*f_k+b\Phi^*g_k)=a\zeta(\Phi^*f_k)+b\zeta(\Phi^*g_k)$ due to the symplecticity of $\zeta$. Fixing $a,b$ and passing to the limit as $k\to\infty$, we see that $\zeta(a\Phi^*f+b\Phi^*g)=a\zeta(\Phi^*f)+b\zeta(\Phi^*g)$, as required. Here we used the fact that quasi-states are Lipschitz relative to the $C^0$-norm, see Remark \ref{rem:Riesz_rep}.

Conversely, if $\Phi_*\zeta$ is linear on $\im\Phi$ whenever $\Phi$ is an involutive map, and we take Poisson-commuting $f,g \in C^\infty(M)$, then the map $\Phi=(f,g)\fc M \to \R^2$ is involutive, and thus $\Phi_*\zeta$ is linear. It then follows that $\zeta$ is linear on $\spn_\R(f,g)=\spn_\R(\Phi^*x,\Phi^*y)$, where $x,y$ are the coordinate functions on $\R^2$.
\end{proof}

\subsection{Quasi-states and topological measures}\label{ss:QS_TM}

Throughout this section, $X$ stands for a compact Hausdorff topological space. We let $\cO(X)$ and $\cC(X)$ denote be the collections of open, respectively closed subsets of $X$, and let $\cA(X) = \cO(X)\cup\cC(X)$. A \tb{topological measure} on $X$ is a function $\tau \fc \cA(X) \to [0,1]$ satisfying
\begin{itemize}
  \item \tb{(normalization):} $\tau(X)=1$;
  \item \tb{(monotonicity):} $\tau(A)\leq\tau(A')$ for $A,A'\in\cA(X)$ with $A\subset A'$;
  \item \tb{(additivity):} $\tau\big(\bigcup_{i=1}^kA_i\big)=\sum_{i=1}^{k}\tau(A_i)$ for pairwise disjoint $A_1,\dots,A_k\in\cA(X)$ with $\bigcup_iA_i\in\cA(X)$;
  \item \tb{(regularity):} for $K\in\cC(X)$, $\tau(K) = \inf\{\tau(U)\,|\,\cO(X)\ni U\supset K\}$.
\end{itemize}
We say that $\tau$ is \tb{simple} if $\tau(\cA(X)) = \{0,1\}$.

Topological measures were introduced by Aarnes \cite{Aarnes_QS_QM}.\footnote{In the early literature they were called quasi-measures and non-subadditive measures.} He also proved \emph{ibid.} a generalization of the Riesz representation theorem (see Remark \ref{rem:Riesz_rep}). Namely, he established a bijection between the collections of topological measures and quasi-states on $X$, which we will now describe. Given a topological measure $\tau$, the corresponding quasi-state is the functional
\begin{equation}\label{eqn:from_TM_to_QS}
  C(X)\ni f\mapsto \min f + \int_{\min f}^{\max f}\tau(\{f\geq s\})\,ds\,.
\end{equation}

Conversely, given a quasi-state $\zeta$, the set function
\begin{equation}\label{eqn:from_QS_to_TM}
  \cC(X)\ni K\mapsto \inf\{\zeta(f)\,|\,C(X)\ni f\geq 1_K\}\,,\quad \text{where }1_K\text{ is the indicator function of }K\,,
\end{equation}
extends to a unique topological measure. We refer to this as the \tb{Aarnes representation theorem}.
\begin{rem}\label{rem:Aarnes_rep_simple}
  Under Aarnes representation, simple quasi-states correspond to simple topological measures \cite{Aarnes_Pure_QS_Extremal_QM} (there simple topological measures are dubbed `extremal quasi-measures').
\end{rem}

In Section \ref{s:intro} we mentioned the importance of $\zeta$-superheavy sets, where $\zeta$ is a quasi-state, see Theorem \ref{thm:sh_Lag_torus_contains_sh_point} and the discussion immediately preceding it. We have the following elementary characterization of these in terms of the representing topological measure.
\begin{lemma}\label{lem:char_sh_sets_tau_is_one}
  Let $X$ be a compact Hausdorff space, let $\zeta$ be a quasi-state on $X$, and assume it is represented by the topological measure $\tau$. Then a compact subset $K\subset X$ is $\zeta$-superheavy if and only if $\tau(K) = 1$.
\end{lemma}
\begin{proof}
  Assume that $K$ is $\zeta$-superheavy, that is for all $c \in \R$ and $f \in C(X)$ with $f|_K\equiv c$ we have $\zeta(f)=c$. Then $\tau(K) = \inf\{\zeta(h)\,|\, h\in C(X)\,, h\geq 1_K\}$. For any $h \geq 1_K$ there exists $f \in C(X)$ such that $1_K \leq f \leq h$ and such that $f|_K\equiv 1$. It follows that $\tau(K) \geq 1$, since $\zeta(f)=1$ for each such $f$. Obviously $\tau(K) \leq 1$, and thus $\tau(K) = 1$. Conversely, if $\tau(K) = 1$ and $f \in C(X)$ is such that $f|_K \equiv c$, let $h \in C(X)$ be such that $h|_K\equiv c$, and such that $h \leq f$ and $h \leq c$. Then
  $$\zeta(f) \geq \zeta(h) = \min h + \int_{\min h}^{\max h}\tau(\{h\geq s\})\,ds\,.$$
  We have $\max h = c$ and for each $s < c$, $K\subset \{h\geq s\}$, therefore $\tau(\{h\geq s\})=\tau(K) = 1$. Thus the integral equals $c-\min h$, therefore $\zeta(f) \geq \zeta(h) = \min h + c-\min h = c$, that is $\zeta(f) \geq c$. Applying this to $-f$, we see that $\zeta(-f) \geq -c$, therefore $\zeta(f) = -\zeta(-f) \leq -(-c) = c$, and the proof is complete.
\end{proof}

In Section \ref{s:intro} we have already seen the notion of pushforwards for quasi-states. It can likewise be defined for topological measures. Namely, if $\Phi \fc X \to Y$ is a continuous map, where $Y$ is another compact Hausdorff space, and $\tau$ is a topological measure on $X$, then the set function $\Phi_*\tau$, given by $\Phi_*\tau(B) = \tau(\Phi^{-1}(B))$ for $B\in\cA(Y)$, is again a topological measure. Moreover, pushforwards commute with the Aarnes representation:
\begin{lemma}\label{lem:pushfwd_commutes_Aarnes_rep}
  If $\zeta$ is a quasi-state on $X$ represented by $\tau$, then $\Phi_*\zeta$ is represented by $\Phi_*\tau$.
\end{lemma}
\begin{proof}
  Let $\eta$ be the quasi-state on $Y$ represented by $\Phi_*\tau$. For $g \in C(Y)$ we have
  $$\eta(g) = \min g + \int_{\min g}^{\max g}\Phi_*\tau(\{g\geq s\})\,ds\,.$$
  Since $\min g = \min (g\circ\Phi)$, $\max g = \max (g\circ\Phi)$, and $\Phi_*\tau(\{g\geq s\}) = \tau(\Phi^{-1}(\{g\geq s\})) = \tau(\{g\circ\Phi\geq s\})$, we see that
  $$\eta(g) = \min(g\circ\Phi) + \int_{\min(g\circ\Phi)}^{\max(g\circ\Phi)} \tau(\{g\circ\Phi\geq s\})\,ds\,,$$
  which by the Aarnes representation theorem equals $\zeta(g\circ\Phi)=\Phi_*\zeta(g)$, see equation \eqref{eqn:from_TM_to_QS}.
\end{proof}

For the proof of our main result we will need the following basic observation, which has already been alluded to in Section \ref{s:intro}.
\begin{lemma}\label{lem:pushfwd_simple_sympl_QS_delta_meas}
  Let $\zeta$ be a simple symplectic quasi-state on a closed symplectic manifold $(M,\omega)$, and let $\Phi\fc M \to \R^k$ be an involutive map. Then $\Phi_*\zeta$ is a delta-measure on $C(\im\Phi)$.
\end{lemma}
\noindent For the proof we need the following elementary
\begin{lemma}\label{lem:simple_meas_delta}
  A regular simple Borel probability measure $\mu$ on $X$ is a delta-measure.
\end{lemma}
\begin{proof}
  Otherwise, for each $x \in X$ we have $\mu(\{x\}) = 0$, which by regularity implies the existence of an open neighborhood $U_x$ of $x$ with $\mu(U_x) = 0$. Since the $U_x$ cover $X$, and since $X$ is compact, there is an open subcover $(U_{x_i})_{i=1}^r$, which implies that $1=\mu(X) \leq \sum_{i=1}^{r}\mu(U_{x_i}) = 0$, which is a contradiction.
\end{proof}
\begin{proof}[Proof of Lemma \ref{lem:pushfwd_simple_sympl_QS_delta_meas}]
  Since $\zeta$ is symplectic and $\Phi$ is involutive, by Lemma \ref{lem:char_sympl_QS_invol_maps}, the quasi-state $\Phi_*\zeta$ on $\im\Phi$ is linear, therefore it is represented by a measure, see Remark \ref{rem:Riesz_rep}. On the other hand if $\tau$ is the topological measure representing $\zeta$, then by Lemma \ref{lem:pushfwd_commutes_Aarnes_rep}, $\Phi_*\zeta$ is represented by the topological measure $\Phi_*\tau$, which, as we have just shown, is in fact a measure. By assumption, $\zeta$, and thus $\Phi_*\zeta$, is simple, and thus so are $\tau$ and $\Phi_*\tau$, see Remark \ref{rem:Aarnes_rep_simple}. Thus $\Phi_*\tau$ is a simple Borel probability measure on the compact space $\im \Phi$, which is regular, since it is a topological measure,  and thus it is a delta-measure by Lemma \ref{lem:simple_meas_delta}. This concludes the proof.
\end{proof}

\begin{defin}\label{def:fiber_finite}
  Let us call a continuous map \tb{fiber-finite} if each one of its fibers has only finitely many connected components. We call such a connected component a \tb{fiber component} of the given map.
\end{defin}

\begin{lemma}\label{lem:fiber_finite_invol_maps_special_fiber_cpt}
  Let $(M,\omega)$ be a closed symplectic manifold, let $\Phi \fc M \to \R^k$ be a fiber-finite involutive map, and let $\zeta$ be a simple symplectic quasi-state on $M$. Then there is a unique fiber component $C_\Phi$ of $\Phi$ satisfying $\tau(C_\Phi)=1$, where $\tau$ is the topological measure representing $\zeta$.
\end{lemma}
\begin{proof}
  By Lemma \ref{lem:pushfwd_simple_sympl_QS_delta_meas}, $\Phi_*\zeta$ is represented by the delta-measure at a point $x_0\in\im\Phi$, therefore $\Phi_*\tau=\delta_{x_0}$ by Lemma \ref{lem:pushfwd_commutes_Aarnes_rep}, whence $\tau(\Phi^{-1}(x_0)) = 1$. Since $\Phi$ is fiber-finite, $\Phi^{-1}(x_0)=C_1\cup\dots\cup C_k$, where the $C_i$ are the distinct connected components of $\Phi^{-1}(x_0)$. Since the $C_i$ are closed, by the additivity of $\tau$ we have $\sum_{i=1}^{k}\tau(C_i)=1$. Since $\tau$ is simple, exactly one of the summands equals $1$, the rest being zero. The corresponding fiber component is the one we are looking for.
\end{proof}
\begin{defin}\label{def:special_fiber_cpt}
  Under the conditions of Lemma \ref{lem:fiber_finite_invol_maps_special_fiber_cpt}, we call $C_\Phi$ the \tb{special fiber component of $\Phi$ with respect to $\tau$}.
\end{defin}

\begin{rem}\label{rem:special_fiber_cpts_intersect}
  Again under the conditions of Lemma \ref{lem:fiber_finite_invol_maps_special_fiber_cpt}, if $\Phi'\fc M \to \R^{k'}$ is another fiber-finite involutive map, then $C_\Phi\cap C_{\Phi'}\neq \varnothing$. Indeed, $\tau(C_\Phi)=\tau(C_{\Phi'})=1$, and thus these two sets must intersect to avoid a contradiction with the additivity and monotonicity of $\tau$.
\end{rem}

We have already defined solid sets in Section \ref{s:intro}. Let us discuss these and their role in proving Theorem \ref{thm:Aarnes_constr} and Corollary \ref{cor:Aarnes_delta}. Solid sets were introduced by Aarnes in \cite{Aarnes_Constr_NSA}, where he also established the fundamental result relating them to topological measures. To formulate it, assume in addition that $X$ is connected and locally connected. \emph{Ibid.}, Aarnes proved that a topological measure is uniquely determined by its values on solid sets in $\cA(X)$, that its restriction to the collection of all such sets is a solid set function, see \cite[Definition 2.3]{Aarnes_Constr_NSA}, and finally that any solid set function uniquely extends to a topological measure.

\begin{exam}\label{ex:solid_sets}
  In the following examples, $X$ is a fixed closed connected manifold.
  \begin{enumerate}
    \item Let $P,P'\subset \C^n$ be polydisks and assume that $\ol{P'}\subset P$. Let $\iota \fc P \hookrightarrow X$ be a smooth embedding, where we assume that $\dim X = \dim P = \dim P'= 2n$. Then $\iota(P')$ is solid. Clearly $\iota(P')$ is connected. Consider the Mayer--Vietoris sequence of the pair $(\ol{\iota(P')},X\setminus\iota(P'))$:
        $$\dots \to H_0(\partial\iota(P')) \to H_0(\ol{\iota(P')}) \oplus H_0(X\setminus\iota(P')) \to H_0(X)\to0\,.$$
        Since $\partial\iota(P')\cong\partial P'$, $X$, and $\ol{\iota(P')}\cong\ol{P'}$ are connected, we see from the exactness that $H_0(X\setminus\iota(P'))$ must have rank at most $1$. Since clearly $\iota(P')\subsetneq X$, it follows that $X\setminus\iota(P')\neq\varnothing$, and thus this rank is at least one. Consequently, $X\setminus \iota(P')$ is connected, and thus $\iota(P')$ is solid as claimed.
    \item If $Q\subset X$ is a closed manifold of codimension at least two, then it does not disconnect $X$. It follows that if $Q$ is itself connected, then it is solid. In particular, if $X$ is symplectic and has dimension at least $4$, and $Q$ is a connected Lagrangian in $X$, then it is solid.
    \item Generalizing item (i), if $U\subset X$ is a connected open set such that $\partial U$ is connected and regular enough to invoke the Mayer--Vietoris sequence for the pair $(\ol U, X\setminus U)$, then $X\setminus U$ is connected. This is the case, for instance, if $\ol U$ is a manifold with corners and connected boundary, such that $U$ is both its topological and manifold interior. We will use the special case when $U$ is a tubular neighborhood of a smooth torus $T\subset X$ such that there exists a diffeomorphism of a neighborhood of $\ol U$ mapping $\ol U$ onto $T \times [-1,1]^n \subset T \times \R^n$, where $n\geq 2$.
  \end{enumerate}
\end{exam}

\paragraph*{Historical remarks.} As we have mentioned in the introduction, Aarnes introduced a notion of genus for $X$ \cite{Aarnes_Constr_NSA}. Knudsen showed in \cite{Knudsen_Topology_constr_extreme_TM} that if $X$ is a finite CW complex, then it has genus zero if and only if $H^1(X;\Z) = 0$. On spaces of genus zero a topological measure can be defined in a very simple manner starting from regular Borel probability measures. To describe it, let us define the \tb{split-spectrum} of such a measure $\mu$ as
$$\{s \in (0,1)\,|\,\exists \text{ solid }C,C' \in \cC(X):C\cap C'=\varnothing\,,\mu(C)=s\,,\mu(C')=1-s\}\,.$$
For a regular Borel probability measure $\mu$ with empty split-spectrum, Aarnes defined in \cite[Section 6, p233]{Aarnes_Constr_NSA} a unique topological measure $\tau$ with the property that for solid $U\in \cO(X)$ we have $\tau(U) = 0$ if $\mu(U)\leq\frac 12$ and $\tau(U)=1$ otherwise (one needs to substitute $n=0$ into the construction appearing in that paper). Later, Butler showed that the same construction works if one only assumes that $\frac12$ does not belong to the split-spectrum of $\mu$ \cite[Theorem 3.5]{Butler_Q_fcns}. \emph{This is the Aarnes construction for topological measures}. The notion of \emph{q-functions} is instrumental for it, as well as for certain generalizations thereof, and it seems that it has first appeared in \cite{Aarnes_Rustad_Prob_QM_new_interpret}. There, the authors present a more general construction of topological measures using q-functions and measures with empty split-spectrum, see Proposition 2.2 \emph{ibid}. Butler's result which we have just mentioned generalizes this as well.

\begin{notation}\label{not:P_0_X}
  For a compact Hausdorff space $X$, we denote by $\cP_0(X)$ be the collection of regular Borel probability measures on $X$ whose split-spectrum does not contain $\frac12$.
\end{notation}

Let us see how this allows us to prove Theorem \ref{thm:Aarnes_constr}. Before doing so, we make the following

\begin{rem}\label{rem:Borel_meas_cpt_metric_space_regular}
  On a compact metrizable space, any finite Borel measure is regular. This follows from \cite[Chapter 6, Proposition 1.3]{Stein_Shakarchi_Real_analysis_3}.
\end{rem}
\begin{proof}[Proof of Theorem \ref{thm:Aarnes_constr}]
  We are given a connected finite CW complex $X$ with $H^1(X;\Z) = 0$, therefore the Aarnes construction for topological measures applies to $\mu$ and yields a topological measure $\tau$. Let us prove that the quasi-state it represents satisfies the asserted property that $\zeta$ vanishes on functions with support in any solid open set $U$ with $\mu(U) \leq \frac12$. Indeed, let $U$ be any such open set, which we may without loss of generality assume nonempty, and consider a nonnegative function $f$ supported in $U$. Since $\mu(X\setminus U) \geq \frac 12$, we have $X\setminus U\neq \varnothing$, and thus $\min f = 0$. Also for any $s > 0$ we have $\{f\geq s\} \subset U$, therefore $\tau(\{f\geq s\}) = 0$, and we conclude from equation \eqref{eqn:from_TM_to_QS} that $\zeta(f) = 0$. Similarly, if $f$ is supported in $U$ and is nonpositive, then $\zeta(f) = -\zeta(-f) = 0$ by the case we have just considered. Since any function $f$ with support in $U$ is bounded between two functions $f_\pm$ supported in $U$ and such that $\pm f_\pm\geq 0$, by the monotonicity of $\zeta$ we conclude that $\zeta(f) = 0$.

  Let us now show uniqueness. Let $\zeta$ be any quasi-state on $X$ vanishing on functions supported in any open solid set $U$ with $\mu(U) \leq \frac 12$, and let $\tau$ be the topological measure representing it. Let $U$ be such an open set. For any compact $K\subset U$ there exists a continuous function $f$ supported in $U$ with $f\geq 1_K$, which implies that $\tau(K) = 0$ by equation \eqref{eqn:from_QS_to_TM}. From the regularity of $\tau$ it follows that $\tau(U) = 0$. Let now $U$ be a solid open set with $\mu(U) > \frac12$. Since $\mu(X\setminus U) < \frac12$, and since $\mu$ is regular, thanks to Remark \ref{rem:Borel_meas_cpt_metric_space_regular}, there exists an open set $V$ containing $X\setminus U$ and such that $\mu(V) < \frac12$. Since $X\setminus U$ is solid, it follows from \cite[Lemma 3.3]{Aarnes_Constr_NSA} that there is an open solid set $W$ such that $X\setminus U\subset W\subset V$. In particular $\mu(W)\leq\mu(V)<\frac12$, and thus $\tau(W)=0$, as we have just proved. It follows that $\tau(U) \geq \tau(X\setminus W) = 1 -\tau(W) = 1$. Thus $\tau$ is the topological measure given by the Aarnes construction for topological measures, and the uniqueness is established.
\end{proof}

\begin{proof}[Proof of Corollary \ref{cor:Aarnes_delta}]
  Let $\tau$ be the topological measure representing the given Aarnes quasi-state $\zeta$, and let us show that $\tau(\{x_0\}) = 1$. Since a point does not disconnect a closed connected manifold, it follows that $\{x_0\}$ is solid. It follows that $X\setminus\{x_0\}$ is a solid open set of $\mu$-measure $\leq\frac 12$, therefore $\tau(X\setminus\{x_0\}) = 0$ and $\tau(\{x_0\})=1$. Thus $\{x_0\}$ is $\zeta$-superheavy by Lemma \ref{lem:char_sh_sets_tau_is_one}, which means that $\zeta$ is indeed the evaluation at $x_0$.
\end{proof}

We close this section with the following elementary observation, which is crucial for the proof of our main result, and which will be used in Section \ref{s:pf_main_result}.
\begin{lemma}\label{lem:solid_cpt_sh_meas_more_half}
  Let $X$ be a finite CW complex with $H^1(X;\Z)=0$, let $\mu\in\cP_0(X)$, and let $\tau$ be the topological measure obtained from $\mu$ by the Aarnes construction. If $K\subset X$ is a solid compact set, then $\tau(K) = 1$ if and only if $\mu(K)\geq\frac12$.
\end{lemma}
\begin{proof}
  Observe that $X\setminus K$ is an open solid set. From the Aarnes construction it follows that $\mu(X\setminus K)\leq \frac12$ if and only if $\tau(X\setminus K)=0$. The claim follows by passing to the complements.
\end{proof}

\section{Constructing the embedding}\label{s:constructing_emb}

This section is dedicated to the proof of Theorem \ref{thm:main_embedding}, which states: Given a closed connected symplectic manifold $(M,\omega)$ and a Borel probability measure $\mu$ on it, for each $\epsilon > 0$ there exists a polydisk $P$ and a symplectic embedding $\iota \fc P \hookrightarrow M$ with $\mu(M\setminus\iota(P))<\epsilon$.

The proof generally follows \cite[Section 6.1]{Schlenk_Embeddings_problems} with appropriate adjustments to the case of a general measure. It is split into three major steps. In the first step, covered in Section \ref{ss:covering_by_cubes}, we construct a symplectic embedding of a finite disjoint union of congruent cubes into $M$, whose image has $\mu$-measure $>1-\epsilon$. In the second step, appearing in Section \ref{ss:embedding_polydisks_cubes}, we construct symplectic embeddings of polydisks into cubes. In the third and final step, detailed in Section \ref{ss:covering_most_meas_polydisk}, we combine the embeddings of the first two steps together with an auxiliary construction in order to construct the desired embedding of a polydisk into $M$.

\subsection{Embedding
cubes}\label{ss:covering_by_cubes}

The goal in this section is to prove
\begin{thm}\label{thm:cover_most_meas_cubes}
  Let $(M,\omega)$ be a closed symplectic manifold of dimension $2n$ and let $\mu$ be a Borel probability measure on $M$. Then for each $\epsilon > 0$ there exists a finite disjoint union $W$ of congruent open cubes in $\R^{2n}$, and a symplectic embedding $\iota \fc W \hookrightarrow M$ such that $\mu\big(M\setminus\iota(W)\big) < \epsilon$.
\end{thm}
\noindent The idea is to first cut out of $M$ a closed set of zero measure, such that its complement decomposes into a finite number of open sets which are small enough to be symplectomorphic to open subsets of $\R^{2n}$. Then inside each such open subset we construct a collection of carefully chosen pairwise disjoint symplectically embedded open cubes, such that their total measure is $>1-\epsilon$. We refer the reader to \cite[Section 6.1]{Schlenk_Embeddings_problems} for a proof in case $\mu$ is the normalized Lebesgue measure corresponding to the volume form $\omega^{\wedge\frac12\dim M}$.

\begin{figure}
  \centering
  \includegraphics[width=8cm]{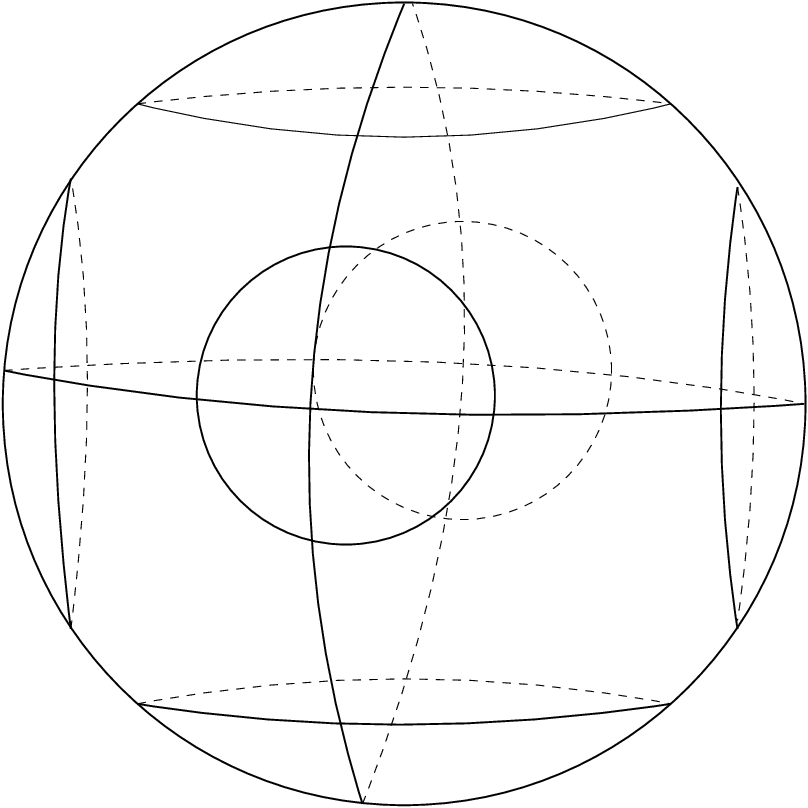}
  \caption{The decomposition $M=C\sqcup\bigsqcup_iU_i$, where $C$ is the union of the curves, while the $U_i$ are the components of the complement.}\label{fig:lattice_sphere}
  \end{figure}

A \tb{Darboux ball} in a symplectic manifold is an open subset symplectomorphic to an open Euclidean ball.
\begin{prop}\label{prop:decomp_sympl_mfd_C_U_i}
  Let $(M,\omega)$ be a closed symplectic manifold and let $\mu$ be a Borel probability measure on $M$. Then there is a decomposition $M = C \sqcup \bigsqcup_{i=1}^k U_i$, where $C\subset M$ is closed, $U_1,\dots,U_k\subset M$ are open, $\mu(C) = 0$, and each $U_i$ is contained in a Darboux ball.
\end{prop}
\begin{proof}
  Fix a finite cover $\cV$ of $M$ by Darboux balls. Without loss of generality we can assume that $M$ is given as a submanifold in some $\R^N$. We equip $M$ with the metric induced from the $\ell^\infty$ metric $d$. Let $\delta$ be a Lebesgue number for $\cV$ relative to this induced metric, that is for each $x \in M$, $B_x(\delta)\cap M$ is contained in one of the members of $\cV$.

  Let $\nu$ be the Borel probability measure induced on $\R^N$ from $\mu$ by the inclusion $M\hookrightarrow \R^N$. By Lemma \ref{lem:Sigma_y_a_at_most_countable}, there is $y \in \R^N$ such that $\nu(\Sigma(y,\delta)) = 0$. Let $\cU = \{U\cap M\,|\,U \in \cU(y,\delta)\,, U\cap M \neq \varnothing\}$. Let $U_1,\dots,U_k$ be the distinct elements of $\cU$, and let $C = M\setminus\bigcup_{i=1}^kU_i = M \cap \Sigma(y,a)$. We claim that $M=C\cup\bigcup_iU_i$ is the required decomposition of $M$. Figure \ref{fig:lattice_sphere} illustrates the decomposition.

  Clearly $C$ is closed, and by construction we have $\mu(C) = \nu(\Sigma(y,\delta))=0$. Each $U_i$ is contained in an open cube with side length $\delta$, which implies that its $d$-diameter is $\leq \delta$, in particular for $x \in U_i$ we have $U_i\subset B_x(\delta)\cap M$, and by assumption this latter set is contained in one of the members of $\cV$.
\end{proof}

We can now present
\begin{proof}[Proof of Theorem \ref{thm:cover_most_meas_cubes}]
  Let $M=C\sqcup\bigsqcup_{i=1}^kU_i$ be a decomposition as in Proposition \ref{prop:decomp_sympl_mfd_C_U_i}. Let $B_i\subset M$ be a Darboux ball containing $U_i$, and let $j_i \fc B_i\hookrightarrow \R^{2n}$ be a symplectic embedding whose image is a Euclidean ball. Applying appropriate shifts, we may assume that the sets $j_i(B_i)$ are pairwise disjoint. Let us define a symplectic embedding $j\fc \bigsqcup_{i=1}^k U_i \hookrightarrow \R^{2n}$ by setting $j|_{U_i}:=j_i|_{U_i}$. Put $U = \im j$ and $\nu=j_*(\mu|_{M\setminus C})$.

  Since $U \subset \R^{2n}$ is open and bounded, and $\nu$ is finite, we have
  $$\sup\{\nu(K)\,|\,K\text{ compact}\subset U\} = \nu(U) = \mu(M\setminus C)= 1\quad \text{by Lemma \ref{cor:approx_open_set_cpts}}.$$
  Let $K\subset U$ be a compact set such that $\mu(K) > 1 - \epsilon$. Put $r:=d(K,\R^{2n}\setminus U)$. Since $K$ is compact, we have $r>0$.

  By Lemma \ref{lem:Sigma_y_a_at_most_countable}, there exists $y \in \R^{2n}$ such that $\nu(\Sigma(y,\frac r2)) = 0$. Fix such $y$ and define
  \begin{equation}\label{eqn:defin_D}
    D = \bigcup_{C \in \cC(y,\tfrac r2):\,C\subset U}C\,.
  \end{equation}
  We claim that $K\subset D$. Indeed, let $x \in K$. It follows that $B_x(r) \subset U$. Consider the closed ball $\ol B{}_x(\tfrac r 4) = \prod_{i=1}^{2n}[x_i - \tfrac r 4,x_i+ \tfrac r 4]$. Clearly $\big(y+(\frac r2\Z)^{2n}\big)\cap \ol B{}_x(\tfrac r 4) \neq \varnothing$. Pick $z$ in this set; it follows that $x \in \ol B{}_z(\tfrac r 4)\subset \ol B_x(\frac r 2) \subset B_x(r)\subset U$, which means that $\ol B{}_z(\tfrac r 4) \in \cC(y,\frac r2)$, and, since $\ol B{}_z(\tfrac r 4) \subset U$, we see that $\ol B{}_z(\tfrac r 4)$ is one of the cubes in the union \eqref{eqn:defin_D} and thus $\ol B{}_z(\tfrac r 4)\subset D$, which finally implies $x \in D$.

  From this we deduce that $\nu(D) \geq \nu(K) > 1 - \epsilon$. Put $W=D\setminus \Sigma(y,r)$. Since $\nu(\Sigma(y,r)) = 0$, it follows that $\nu(W) > 1-\epsilon$. Note that $W$ is a finite disjoint union of members of $\cU(y,r)$, that is a finite disjoint union of cubes with side length $r$. The desired symplectic embedding $\iota \fc W\hookrightarrow M$ is then obtained by $\iota:=j^{-1}|_W$. To finish the proof, we note that $\mu(M\setminus\im\iota) = \mu((M\setminus C)\setminus\im\iota) = \nu(U\setminus W) < \epsilon$.
\end{proof}

\subsection{Embedding polydisks into cubes}\label{ss:embedding_polydisks_cubes}

The goal here is to prove a foundational result for symplectically embedding polydisks into cubes. Recall Notation \ref{not:cube} where we defined the symplectic cube $C^{2n}(a):=(0,a)^{2n}\subset \R^{2n}$ with side $a>0$. Then, denoting by $\Sigma:=\Sigma\big((\tfrac 12,\dots,\tfrac 12),1\big)$ the union of the integer coordinate planes in $\R^{2n}$, we have:
\begin{thm}\label{thm:polydisk_into_cube}
  Let $k\geq 5$ be an odd integer. Then there exists $\alpha > 0$ and a symplectic embedding $\iota\fc (0,\alpha) \times (0,1)^{2n-1} \hookrightarrow C^{2n}(k)$ such that
  \begin{enumerate}
    \item The image of $\iota$ contains $(2,k-2)\times(0,k)^{2n-1} \setminus \Sigma$;
    \item The restriction of $\iota$ to $(0,1)\times(0,1)^{2n-1}$ is the identity;
    \item The restriction of $\iota$ to $(\alpha-1,\alpha)\times(0,1)^{2n-1}$ is the shift $T_{(k-\alpha,\,k-1,\,\dots,\,k-1)}$.
  \end{enumerate}
\end{thm}

\begin{rem}\label{rem:box_symplecto_polydisk}
  Note that $(0,\alpha)\times(0,1)^{2n-1}$ is symplectomorphic to a polydisk. Indeed, $(0,\alpha)\times(0,1)$ is symplectomorphic to $D(\frac{\alpha}{\pi})$, while $(0,1)^2=C^2(1)$ is symplectomorphic to $D(\frac1{\pi})$, therefore $(0,\alpha)\times(0,1)^{2n-1}$ is symplectomorphic to the polydisk $D(\frac{\alpha}{\pi})\times D(\frac1{\pi})^{n-1} = P(\frac{\alpha}{\pi},\frac1{\pi},\dots,\frac1{\pi})$.
\end{rem}

The rest of this subsection is dedicated to the proof of Theorem \ref{thm:polydisk_into_cube}. Before we give the details of the construction, let us present an overview. All the ideas are contained in Schlenk's book \cite{Schlenk_Embeddings_problems}, and in fact there is almost nothing new in our proof, however, since it is somewhat challenging to extract from the book the exact results we need, and since we will use a particular variant of the construction described there, we decided to give a detailed account.

We view $\R^{2n}$ as the product of the ``base'' $\R^2$, whose coordinates will be denoted by $(u,v)$, and the ``fiber'' $\R^{2n-2}$ with coordinates $(x,y)=(x_1,y_1,\dots,x_{n-1},y_{n-1})$. The set $(0,\alpha)\times(0,1)^{2n-1}$ splits as the product of the rectangle $(0,\alpha)\times(0,1)$ in the base and the cube $(0,1)^{2n-2}$ in the fiber. The desired embedding is constructed in three steps, and is obtained as the composition $\iota=\Gamma\circ\Psi\circ(\Theta\times\id_{\R^{2n-2}})$.

The first map, which we refer to as \emph{preparing the base}, is the product $\Theta\times\id_{\R^{2n-2}}$, where $\Theta$ is a certain symplectomorphism from $(0,\alpha)\times(0,1)$ onto an open band consisting of the union of $k^{2n-2}$ rectangles of the form $(ki+2,ki+k-2)\times(0,k)$, $i=0,\dots,k^{2n-2}-1$, and bridges between them.

The second map $\Psi$, called \emph{symplectic lift}, is a certain symplectomorphism such that the projections of all the sets $\Psi(((i-1)k,ik-1)\times(0,k)\times(0,1)^{2n-2})$, $i=1,\dots,k^{2n-2}$, to the fiber $\R^{2n-2}$ are pairwise disjoint. The word `lift' here refers to the shift effected by $\Psi$ in the fiber over each rectangle.

The third and final map $\Gamma$ is given by \emph{symplectic folding}, and has the form $\Gamma=\gamma\times\id_{\R^{2n-2}}$, where $\gamma \fc \im\Theta\to (0,k)^2$ is a certain symplectic immersion in the base. Even though $\Gamma$ is not an embedding, the composition $\iota$ is. Let us now describe all of this precisely.

\subsubsection{Preparing the base}\label{sss:prep_base}

\begin{figure}
  \centering
  \includegraphics[width=12cm]{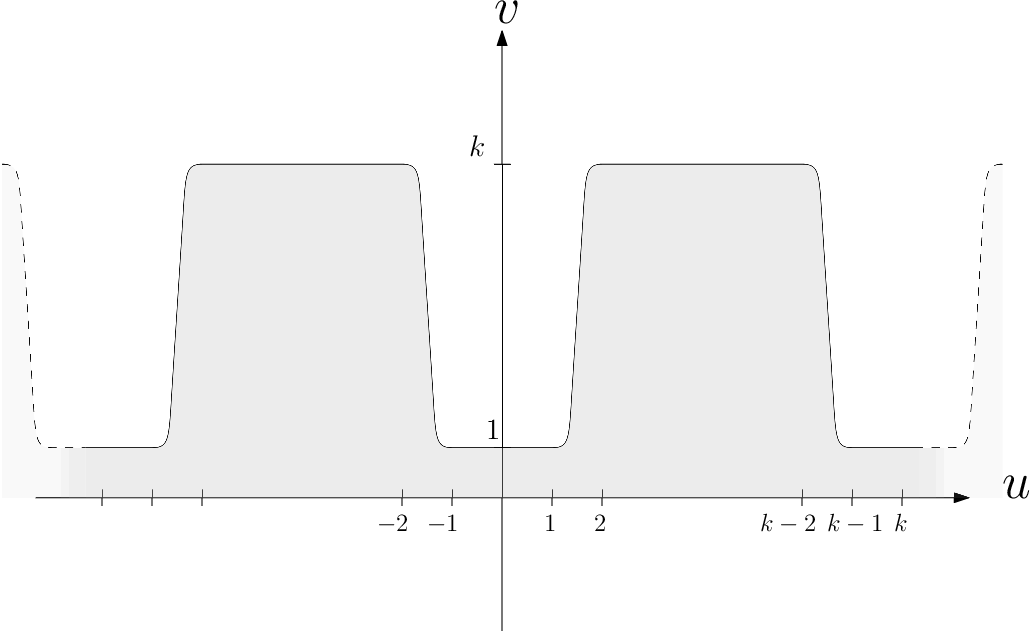}
  \caption{The graph of the function $h$ (Section \ref{sss:prep_base})}\label{fig:graph_h_prep_base}
\end{figure}

Let us fix a smooth $k$-periodic even function $h \fc \R \to [0,k]$, such that $h|_{[0,1]}\equiv 1$, $h|_{[2,k-2]}\equiv k$, and such that on $[1,2]$ it is strictly increasing, see Figure \ref{fig:graph_h_prep_base}. Let $A=\int_{0}^{k}h(t)\,dt$ and $\alpha=k^{2n-2}A$. Let $w \fc \R \to (0,\infty)$ be defined by $w=1/f'$, where $f \fc \R \to \R$ is given by the implicit relation $f^{-1}(x) = \int_0^xh(t)\,dt$ for all $x \in \R$. Thanks to Remark \ref{rem:smears}, the symplectic smear $\Theta_w$ maps the strip $\R\times(0,1)$ symplectomorphically onto the open band bounded by $\R\times\{0\}$ from below and by the graph of $h$ from above. Moreover, since for each $i \in \Z$, $h|_{[ki-1,ki+1]}\equiv 1$, the restriction of $\Theta_w$ to $[iA-1,iA+1]\times \R$ coincides with the shift operator $T_{(i(k-A),0)}$. It follows that $\Theta_w\big((0,\alpha)\times(0,1)\big)$ is the portion contained in $(0,k^{2n-1})\times \R$ of the open band bounded by $\R\times\{0\}$ from below and the graph of $h$ from above.

Let $g \fc \R \to \R$ be the following smooth function: $g(u)\equiv 0$ for $u \leq k^{2n-1}-2$, $g(u)\equiv k-1$ for $u\geq k^{2n-1}$, and $g(u)=k-h(u)$ otherwise. We then define the desired symplectomorphism in the base $\Theta \fc (0,\alpha)\times(0,1) \to \R^2$ as the composition $\Theta=S_g\circ\Theta_w|_{(0,\alpha)\times(0,1)}$, where $S_g$ is the symplectic shear corresponding to the function $g$, see Example \ref{ex:symplectomorphisms} for its definition.

\begin{figure}
  \centering
  \includegraphics[width=12cm]{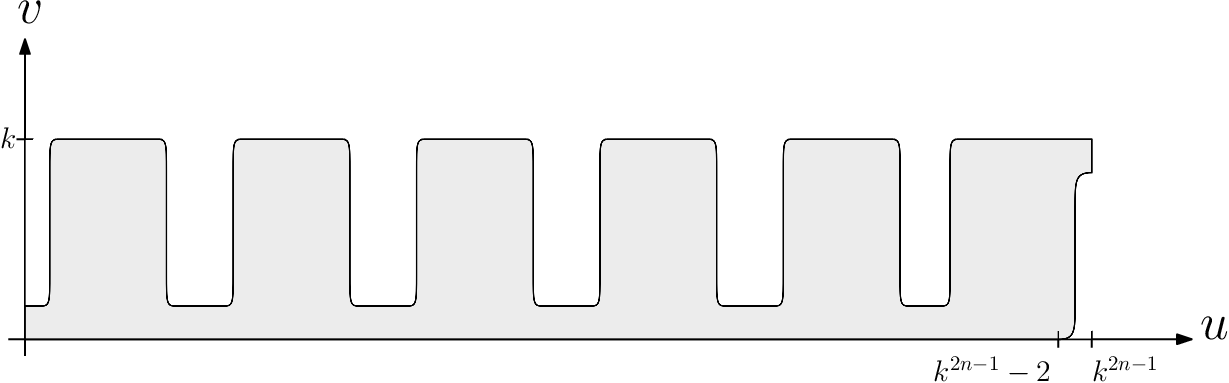}
  \caption{The image of $\Theta$}\label{fig:image_Theta}
\end{figure}

Figure \ref{fig:image_Theta} illustrates the image of $\Theta$, which differs from $\Theta_w\big((0,\alpha)\times(0,1)\big)$ at the last `tail.' More precisely, $\im\Theta$ is the open band bounded by $u=0$ on the left, $u=k^{2n-1}$ on the right, the graph of $g$ on the bottom, and the graph of $\wt h$ on the top, where $\wt h \fc [0,k^{2n-1}] \to \R$ is defined by $\wt h(u)=h(u)$ for $u \leq k^{2n-1}-2$ and by $\wt h\equiv k$ for $u\in[k^{2n-1}-2,k^{2n-1}]$.

\subsubsection{The lift}\label{sss:lift}

In preparation for the last step of the construction, which is folding, we will apply a lifting symplectomorphism. Intuitively speaking, during folding we will fill out the cube $(0,k)^{2n-2}$ in the fiber coordinates by unit cubes, with the exception of some subsets of the union of the integer coordinate planes. In order for the resulting embedding to satisfy the conditions of Theorem \ref{thm:polydisk_into_cube}, this filling process has to start at the cube $(0,1)^{2n-2}$ and terminate at the cube $(k-1,k)^{2n-2}$, where we refer to items (ii, iii) in the formulation of the theorem.

Formally, let us consider the bijection between $\Z^{2n-2}$ and the set of open integral unit cubes in $\R^{2n-2}$, whereby $w \in \Z^{2n-2}$ corresponds to the cube $C_w = \prod_{j=1}^{2n-2}(w_j,w_j+1)$. Under this bijection, the unit cubes contained in $C^{2n-2}(k)=(0,k)^{2n-2}$ correspond to the set $\{0,\dots,k-1\}^{2n-2}\subset\Z^{2n-2}$.

\begin{figure}
  \centering
  \includegraphics[width=8cm]{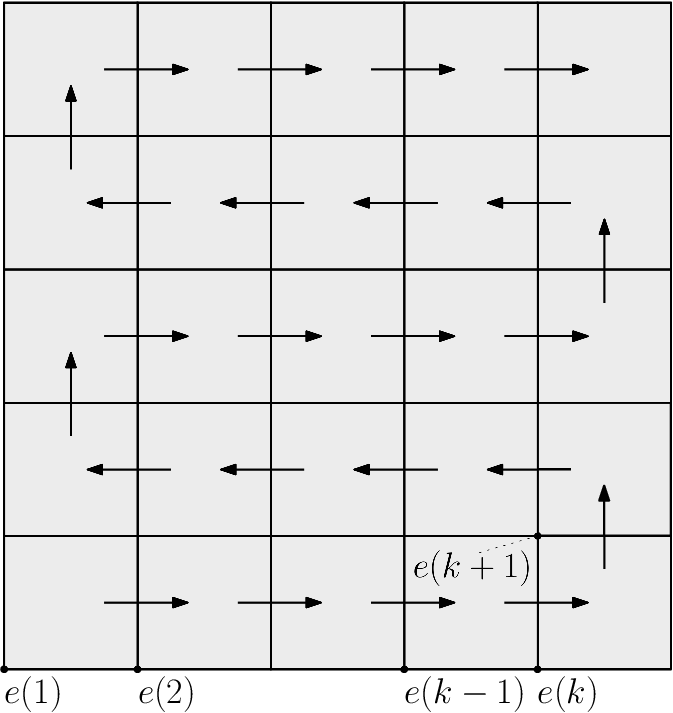}
  \caption{One choice of the enumeration $e$}\label{fig:enum_cubes}
\end{figure}

Now we need an enumeration of $\{0,\dots,k-1\}^{2n-2}$, that is a bijection
$$e \fc \{1,\dots,k^{2n-2}\} \to \{0,\dots,k-1\}^{2n-2}\,,$$
such that $e(1) = (0,\dots,0)$, $e(k^{2n-2}) = (k-1,\dots,k-1)$, and such that for $i \in \{1,\dots,k^{2n-2}-1\}$, the cubes $C_{e(i)}$ and $C_{e(i+1)}$ are adjacent. It is easy to show that such a function exists, for instance we can enumerate the cubes in a `zig-zag' pattern. The assumption that $k$ is odd is crucial here. See Figure \ref{fig:enum_cubes} for an illustration of a possible choice of $e$.

The adjacency condition can be equivalently expressed as follows: for each such $i$, $e(i+1)-e(i)=\pm\partial_{x_j}$ or $\pm\partial_{y_j}$ for some $j=1,\dots,n-1$ depending on $i$. For the lift we also need to record the range of the dual coordinate in the cube $C_{e(i)}$, where the coordinates $x_j$ and $y_j$ are dual to each other. Consequently, we define the following four families of subsets of $\{1,\dots,k^{2n-2}-1\}$, where $j=1,\dots,n-1$ and $\ell=0,\dots,k-1$:
$$X^{\pm}_{j,\ell} = \{i\in\{1,\dots,k^{2n-2}-1\}\,|\,e(i+1) - e(i) = \pm\partial_{x_j}\text{ and }y_j\in(\ell,\ell+1)\text{ in }C_{e(i)}\}\,,$$
$$Y^{\pm}_{j,\ell} = \{i\in\{1,\dots,k^{2n-2}-1\}\,|\,e(i+1) - e(i) = \pm\partial_{y_j}\text{ and }x_j\in(\ell,\ell+1)\text{ in }C_{e(i)}\}\,.$$
By our assumptions on $e$, $\{1,\dots,k^{2n-2}-1\}$ is the disjoint union of these sets as $j,\ell$ run over their ranges.

We will now define the Hamiltonians which will be used in the lifting construction. Pick a smooth function $\rho \fc \R\to [0,1]$ such that $\rho$ is nondecreasing, $\rho|_{[-\frac 14,\infty)}\equiv 1$, $\rho_{(-\infty,-\frac 34]}\equiv 0$, and such that $\rho'<    3$ everywhere.

\begin{figure}
  \centering
  \includegraphics[width=15cm]{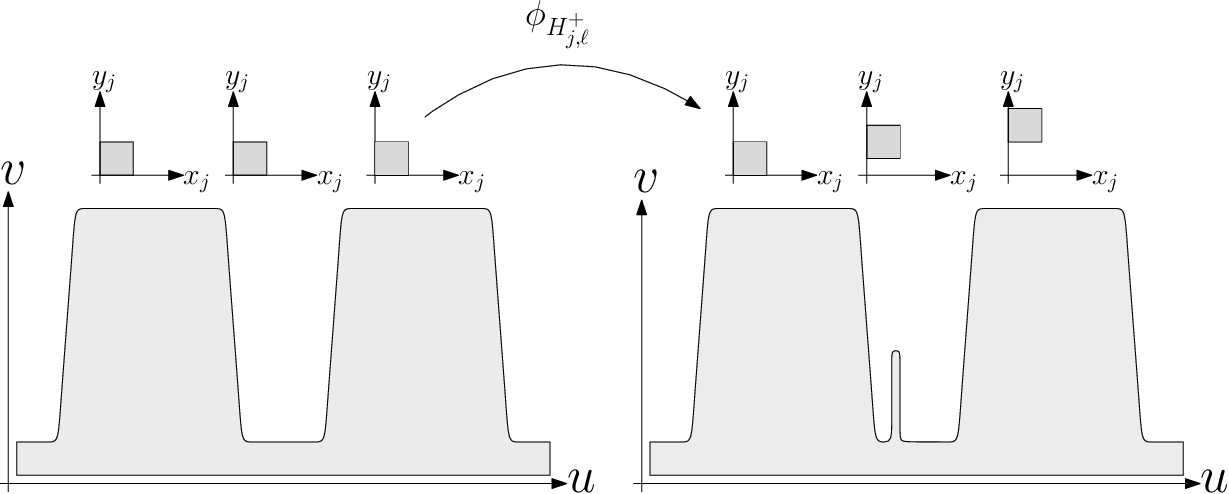}
  \caption{The action of the lifting maps $\phi_{H^+_{j,\ell}}$. The little graphs in the top part are the cross-sections, in the coordinate pair $x_j,y_j$, of the set $\im\Theta \times (0,1)^{2n-2}$ before and after the action}\label{fig:lifting_maps}
\end{figure}

The Hamiltonians come in four families, corresponding to the sets $X^\pm_{j,\ell},Y^\pm_{j,\ell}$. Namely, for $j=1,\dots,n-1$, and $\ell=0,\dots,k-1$, we define:
$$G^+_{j,\ell}(u,v;x,y) = \rho(u)(\ell+1-y_j)\,,\quad G^-_{j,\ell}(u,v;x,y)=\rho(u)(y_j-\ell)\,,$$
$$H^+_{j,\ell}(u,v;x,y) = \rho(u)(x_j-\ell)\,,\quad H^-_{j,\ell}(u,v;x,y)=\rho(u)(\ell+1-x_j)\,.$$
See Example \ref{ex:symplectomorphisms}, `Symplectic lifts' in Section \ref{ss:sympl_geom}. It follows that the time-$1$ maps of these Hamiltonians are:
$$\phi_{G^+_{j,\ell}}(u,v;x,y) = (u,v;x,y)+\rho'(u)(\ell+1-y_j)\partial_v + \rho(u)\partial_{x_j}\,,$$
$$\phi_{G^-_{j,\ell}}(u,v;x,y) = (u,v;x,y)+\rho'(u)(y_j-\ell)\partial_v - \rho(u)\partial_{x_j}\,,$$
$$\phi_{H^+_{j,\ell}}(u,v;x,y) = (u,v;x,y)+\rho'(u)(x_j-\ell)\partial_v + \rho(u)\partial_{y_j}\,,$$
$$\phi_{H^-_{j,\ell}}(u,v;x,y) = (u,v;x,y)+\rho'(u)(\ell+1-x_j)\partial_v - \rho(u)\partial_{y_j}\,.$$
That is, we see that these maps are the identity for $u\leq -\frac34$, and that for $u\geq-\frac14$ they are the shift by $\pm\partial_{x_j},\pm\partial_{y_j}$, according to the type. Moreover, when restricted to the slice $\{\ell< x_j <\ell+1\}$, the time-$1$ maps of $H^\pm_{j,\ell}$ shift $v$ by a quantity contained in $[0,\rho'(u)]$, while the same is true of the time-$1$ maps of $G^\pm_{j,\ell}$ restricted to the slice $\{\ell<y_j<\ell+1\}$. Note as well that a nontrivial shift in $v$ only happens for $u \in (-\frac34,-\frac14)$, since $\rho'$ vanishes outside of this interval. See Figure \ref{fig:lifting_maps} for a picture.

\begin{figure}
  \centering
  \includegraphics[width=8cm]{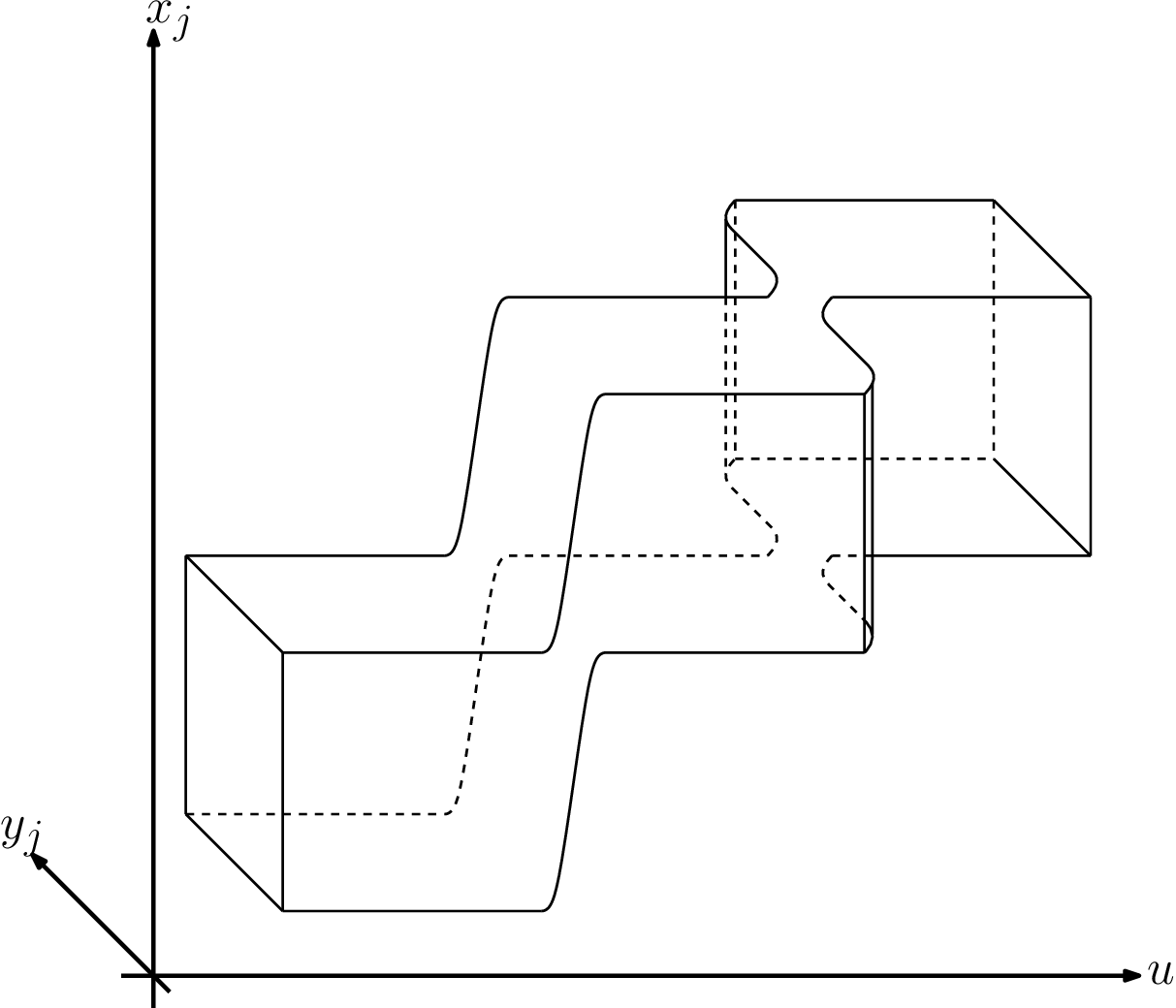}
  \caption{The action of the lifting map $\Psi$ in the coordinates $u,x_j,y_j$}\label{fig:Psi_projected}
\end{figure}

For $i\in\{1,\dots,k^{2n-2}-1\}$, we define a symplectomorphism $\Psi_i$ of $\R^{2n}$ as follows:
\begin{itemize}
  \item If $i \in X^+_{j,\ell}$, then $\Psi_i = T_{ik\partial_u}\circ\phi_{G^+_{j,\ell}}\circ T_{-ik\partial_u}$;
  \item If $i \in X^-_{j,\ell}$, then $\Psi_i = T_{ik\partial_u}\circ\phi_{G^-_{j,\ell}}\circ T_{-ik\partial_u}$;
  \item If $i \in Y^+_{j,\ell}$, then $\Psi_i = T_{ik\partial_u}\circ\phi_{H^+_{j,\ell}}\circ T_{-ik\partial_u}$;
  \item If $i \in Y^-_{j,\ell}$, then $\Psi_i = T_{ik\partial_u}\circ\phi_{H^-_{j,\ell}}\circ T_{-ik\partial_u}$.
\end{itemize}
The total lifting symplectomorphism is the composition
$$\Psi = \Psi_{k^{2n-2}-1}\circ\Psi_{k^{2n-2}-2}\circ\dots\circ\Psi_2\circ\Psi_1\,.$$
It can be checked that for each $i\in\{1,\dots,k^{2n-2}\}$, the projection to the fiber coordinates of the set $\Psi\big(((i-1)k,ik-1)\times(0,k)\times(0,1)^{2n-2}\big)$ is the cube $C_{e(i)}$, which is the chief property of the construction. The projection of the image of $\Psi$ onto the coordinates $u,x_j,y_j$ is depicted in Figure \ref{fig:Psi_projected}.

\subsubsection{The folding}\label{sss:fold}

\begin{figure}
  \centering
  \includegraphics[width=15cm]{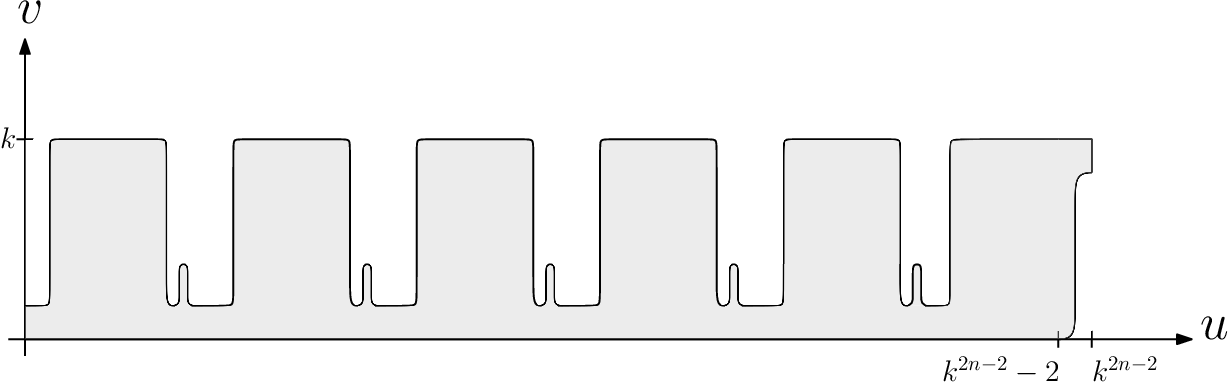}
  \caption{The projection of $\Psi\big(\im\Theta\times(0,1)^{2n-2}\big)$ to $\R^2(u,v)$}\label{fig:projection_Psi_u_v}
\end{figure}

Let $\sigma\fc \R\to\R$ be the $k$-periodic extension of $\rho'|_{[-k,0]}$ to $\R$. It can then be checked that the projection to the base coordinates of the set $\Psi\big(\im\Theta\times(0,1)^{2n-2}\big)$ is the open band bounded by $u=0,u=k^{2n-1}$ on the left and on the right, by the graph of $g$ (see Section \ref{sss:prep_base} for its definition) on the bottom, and on top by the graph of $h+\sigma$ for $u\in(0,k^{2n-1}-2]$ and the line $v=k$ for $u\in[k^{2n-1}-2,k^{2n-1})$, see Figure \ref{fig:projection_Psi_u_v}.

In particular for each $i\in\{1,\dots,k^{2n-2}-1\}$, the intersection of this projection with the rectangle $(ik-1,ik)\times(0,k)$ is contained between $v=0$ and the graph of $1+\sigma$. By our assumption on $\rho$, $1+\sigma<4$, so in this rectangle, there is some space left above this graph, and it will be enough to contain the image of a map we will shortly describe.

For the folding we will only need one additional map, namely $\gamma_+ \fc (-\frac14,\frac14)\times(0,1)\to\R^2$. We claim that there exists such $\gamma_+$ which satisfies the following properties:
\begin{itemize}
  \item $\gamma_+$ is a symplectomorphism onto its image;
  \item $\gamma_+$ is the identity when restricted to $(-\frac14,-\frac18]\times(0,1)$;
  \item the restriction of $\gamma_+$ to $[\frac18,\frac14)\times(0,1)$ is $(u,v) \mapsto (0,k)-(u,v)$;
  \item $\gamma_+([-\frac 18,\frac18]\times(0,1))\subset [-\frac 18,0)\times(0,k)$.
\end{itemize}
\noindent The construction of such a map follows the basic ideas described in Schlenk's book, see \cite[Step 4, p44]{Schlenk_Embeddings_problems}. We will need another map $\gamma_- \fc(-\frac 14, \frac14)\times(0,1)\to\R^2$, defined by $\gamma_-(u,v) = (0,k)-\gamma_+(u,v)$.

\begin{figure}
  \centering
  \includegraphics[width=14cm]{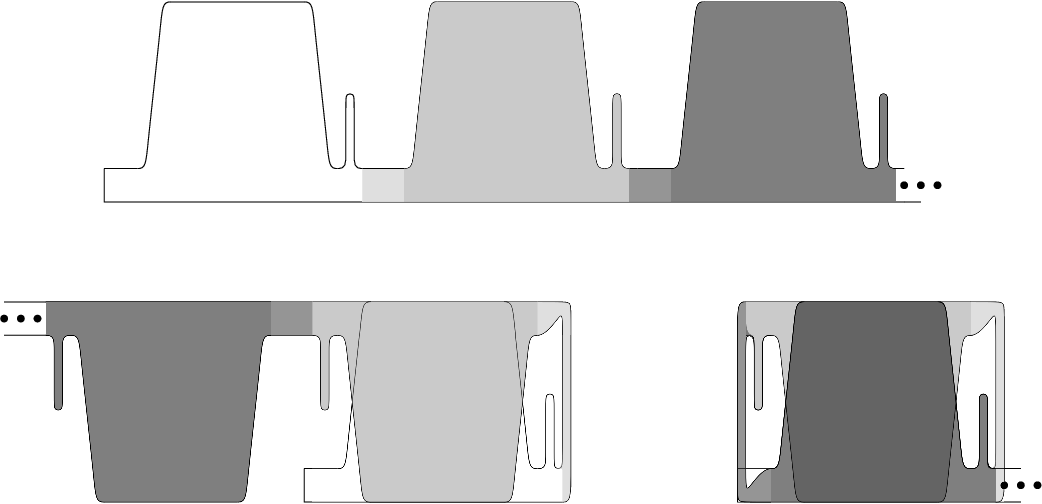}
  \caption{A illustration of the action of $\gamma$. One can think of $\gamma$ as being constructed in stages, where in each stage we take the remaining piece of the band and wrap it around the central square. Top: before the wrapping. Bottom left: after one cycle of wrapping; bottom right: after two cycles of wrapping; here the image of $\gamma$ can already be seen as the union of the shaded areas}\label{fig:gamma}
\end{figure}

Let $\gamma \fc \im\Theta \to \R^2$ be defined as follows:
$$\gamma(u,v)=\begin{cases}(u,v)\,, & u\in(0,k-\frac 18] \\
                           (u-(k^{2n-2}-1)k,v)\,,& u\in[(k^{2n-2}-1)k+\frac 18,k^{2n-2}\cdot k)\\
                           (u-ik,v)&u\in[ik+\frac18,(i+1)k-\frac 18]\,,i\in\{2,\dots,k^{2n-2}-3\}\text{ even}\\
                           ((i+1)k - u,k-v)&u\in[ik+\frac18,(i+1)k-\frac 18]\,,i\in\{1,\dots,k^{2n-2}-2\}\text{ odd}\\
                           (k,0)+\gamma_+(u-ik,v)&u\in(ik-\frac14,ik+\frac 14)\,,i\in\{1,\dots,k^{2n-2}-2\}\text{ odd}\\
                           \gamma_-(u-ik,v)&u\in(ik-\frac14,ik+\frac 14)\,,i\in\{2,\dots,k^{2n-2}-1\}\text{ even}
\end{cases}$$
It is easy to show that $\gamma$ is a symplectic immersion with image $\subset(0,k)^2$. Its action is illustrated in Figure \ref{fig:gamma}.

\begin{figure}
  \centering
  \includegraphics[width=10cm]{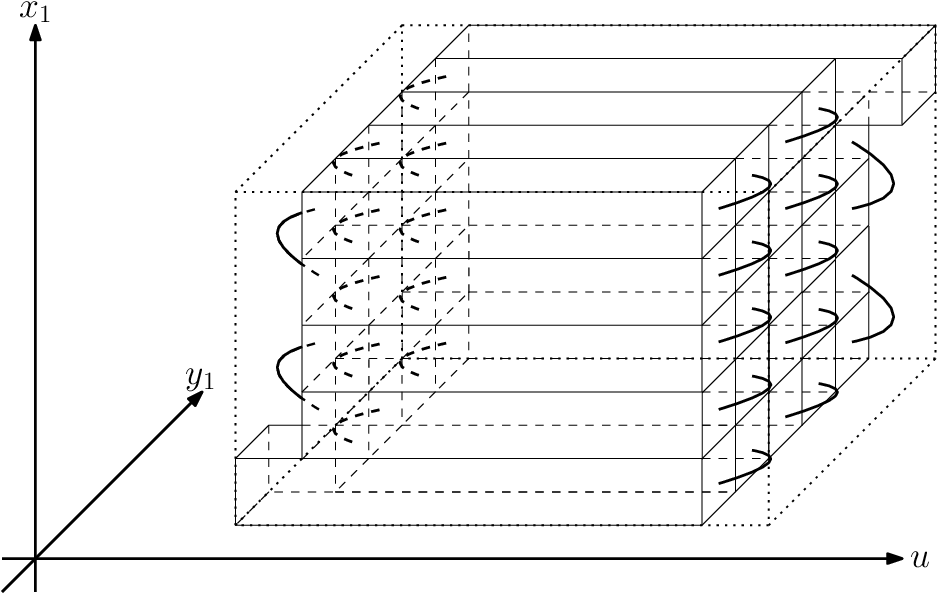}
  \caption{A schematic depiction of the projection of the total embedding map $\iota$ to the coordinates $u,x_1,y_1$}\label{fig:folded_polydisk}
\end{figure}

The folding map $\Gamma$ is defined by $\Gamma=\gamma\times\id_{\R^{2n-2}}$. Note that even though $\gamma$, and therefore $\Gamma$, is not injective, the total composition
$$\iota=\Gamma\circ\Psi\circ(\Theta\times\id)|_{(0,\alpha)\times(0,1)^{2n-1}}$$
is injective, and thus a symplectic embedding, see Remark \ref{rem:injective_sympl}. This is because when $\Gamma$ is applied to different squares in the base, they are already separated in the fiber, which makes the composition injective. Figure \ref{fig:folded_polydisk} illustrates the projection of the image of $\iota$ to the coordinates $u,x_1,y_1$. This concludes the construction of the desired embedding and the proof of Theorem \ref{thm:polydisk_into_cube}.

\subsection{Covering most of the measure by a polydisk}\label{ss:covering_most_meas_polydisk}

Here we will use the constructions of Sections \ref{ss:covering_by_cubes}, \ref{ss:embedding_polydisks_cubes} to prove Theorem \ref{thm:main_embedding}. We keep the notation $(u,v)$ for the base coordinates, and $(x,y)=(x_1,y_1,\dots,x_{n-1},y_{n-1})$ for the fiber coordinates.

The proof is a lengthy construction, split into a number of steps. Let us briefly describe these before passing to details:
  \begin{enumerate}
    \item In \tb{Step 1} we use Theorem \ref{thm:cover_most_meas_cubes} to symplectically embed a finite disjoint union of congruent open cubes in $\R^{2n}$ into $M$ such that the image has measure $>1-\epsilon$. At this point we also pick slightly smaller cubes so that their image still has measure $>1-\epsilon$. This is done to make room for subsequent constructions.
    \item Since the above cubes are located at arbitrary positions in $\R^{2n}$, in \tb{Step 2} we shift them so that they assume the necessary positions for the construction of Step 3.
    \item In \tb{Step 3} we construct thin tubes connecting suitable corners of the smaller cubes, as well as a symplectic embedding of the connected union of the larger cubes and the thin tubes into $M$.
    \item In \tb{Step 4} we choose a suitable thickness for our eventual embedded polydisk, and embed polydisks having this thickness into each cube.
    \item The polydisks constructed in Step 4 do not align in the coordinates $v,x,y$, and in \tb{Step 5} we perform a simple shift to correct this.
    \item In Step 1 we embedded a disjoint union of cubes in $M$, whose image has measure $>1-\epsilon$, and now we have embedded a thin polydisk into each such cube. The issue now is that the image of this polydisk inside each cube \emph{a priori} has unknown measure. Therefore in \tb{Step 6} we apply very small shifts to the polydisks inside each cube so that the result covers most of the measure of the respective cubes.
    \item In \tb{Step 7} we assemble the embeddings constructed in the previous steps to obtain the desired embedding of a polydisk into $M$ covering a set of measure $>1-\epsilon$.
  \end{enumerate}
Let us now describe all of this in details.

\begin{proof}[Proof of Theorem \ref{thm:main_embedding}]Of course, if $\epsilon > 1$, then there is nothing to prove, so in what follows we assume $\epsilon \leq 1$.

  \tb{Step 1.} Theorem \ref{thm:cover_most_meas_cubes} yields some $a>0$ and a symplectic embedding
  $$\iota \fc \bigsqcup_{z \in Z} B_z(a) \hookrightarrow M$$
  with $\mu(\im \iota) > 1-\epsilon$, where $Z\subset \R^{2n}$ is a finite set, and where the cubes are pairwise disjoint. Note that since the empty set has zero measure, to achieve $\mu(\im\iota)>1-\epsilon$ we need at least one cube, meaning $Z$ must be nonempty. Moreover, if $Z$ only has one element, the proof of the theorem stops here. Indeed, since a cube is symplectomorphic to a polydisk, see Example \ref{ex:symplectomorphisms}, the desired embedding is obtained as the composition of $\iota$ with such a symplectomorphism, and its image has measure $>1-\epsilon$, as required. \emph{Henceforth we assume that $Z$ has at least two elements.} By Lemma \ref{cor:approx_open_set_cpts}, there exists $a'' \in (0,a)$ such that
  $$\textstyle\mu \Big(\iota \big(\bigsqcup_{z \in Z}\ol B{}_z(a'')\big)\Big) > 1-\epsilon\,.$$
  We also choose $a' \in (a'',a)$ such that $a' > a-a'$.

  \tb{Step 2.} We now shift all the cubes so that they assume the necessary position for the construction of connecting lines. Let $Z = \{z^{(i)}\}_{i=1}^r$, and note that $r\geq 2$ by our assumption in Step 1. For each $i = 1,\dots,r$ let
  $$w^{(i)}=((4i-3)a',(2i-1)a',\dots,(2i-1)a')\,.$$
  We define a new embedding
  $$\iota' \fc \bigsqcup_{i=1}^r B_{w^{(i)}}(a) \to M\quad\text{via}\quad \iota'|_{B_{w^{(i)}}(a)}=\iota|_{B_{z^{(i)}}(a)}\circ T_{z^{(i)} - w^{(i)}}|_{B_{w^{(i)}}(a)}\,.$$
  Note that the condition $a'>a-a'$ ensures that the cubes $B_{w^{(i)}}(a)$ are pairwise disjoint. See Figure \ref{fig:cubes_R2n} for an illustration, where the first two of the cubes $B_{w^{(i)}}(a')$ are marked with a light shading.

\begin{figure}
  \centering
  \includegraphics[width=13cm]{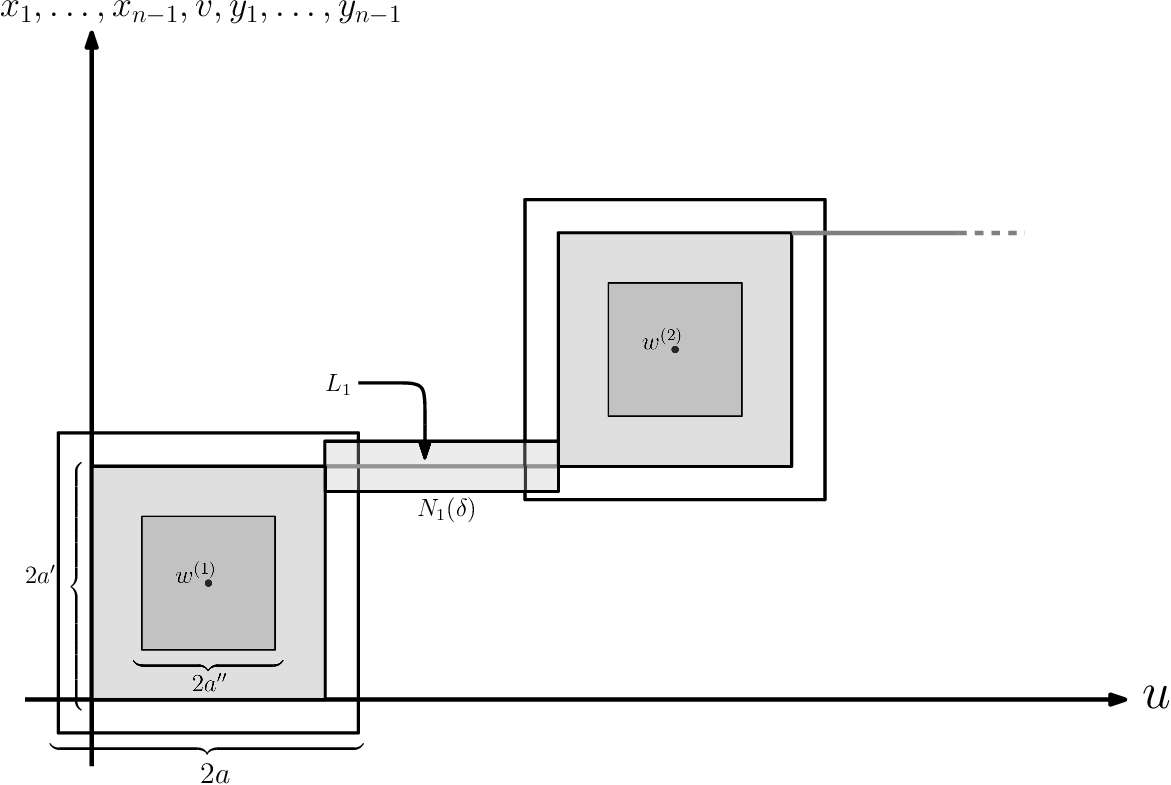}
  \caption{Steps 2 and 3}\label{fig:cubes_R2n}
\end{figure}

  \tb{Step 3.} For $i=1,\dots,r-1$ let $L_i$ be the straight segment connecting
  $$w^{(i)} + (a',\dots,a') = ((4i-2)a',2ia',\dots,2ia')\text{ to } w^{(i+1)}-(a',\dots,a') = (4ia',2ia',\dots,2ia')\,.$$
  For $\delta \in (0,a-a')$ define the corresponding box neighborhood
  $$N_i(\delta) = ((4i-2)a',4ia')\times(2ia'-\delta,2ia'+\delta)^{2n-1}\,.$$
  See Figure \ref{fig:cubes_R2n}, where one can see the first connecting segment $L_1$ and the corresponding box neighborhood $N_1(\delta)$. For $\delta$ small enough we can repeat the construction of Schlenk (see \cite[Proposition 6.1.3]{Schlenk_Embeddings_problems} and the discussion that follows), which yields a symplectic embedding
  $$\iota'' \fc \bigcup_{i=1}^r B_{w^{(i)}}(a) \cup \bigcup_{i=1}^{r-1}N_i(\delta) \hookrightarrow M\,.$$

  \tb{Step 4.} Pick a positive odd integer $k \geq 5$ large enough so that
  \begin{equation}\label{eqn:estimate_eta}
    \eta:=\frac{2a'}{k} < \min\left(\frac{a'-a''}{r-1},\frac{a'-a''}{2},\frac{\delta}{r-1}
    \right)\,.
  \end{equation}
  Note that the second condition is only relevant if $r = 2$.
  Scaling the embedding of Theorem \ref{thm:polydisk_into_cube}, we obtain $\alpha>0$ and a symplectic embedding
  $$\iota_1 \fc (0,\alpha)\times (0,\eta)^{2n-1}\hookrightarrow C^{2n}(k\eta)=(0,2a')^{2n} = B_{w^{(1)}}(a')\,.$$
  Applying translations, we obtain for $i=2,\dots,r$:
  $$\iota_i=T_{w^{(i)}-w^{(1)}}\circ\iota_1 \fc (0,\alpha)\times (0,\eta)^{2n-1}\hookrightarrow B_{w^{(i)}}(a')\,.$$

  \tb{Step 5.} In this step we will use the $\iota_i$ to produce an embedding
  $$\wh\iota \fc (0,\beta) \times (0,\eta)^{2n-1}\hookrightarrow \bigcup_{i=1}^r B_{w^{(i)}}(a) \cup \bigcup_{i=1}^{r-1}N_i(\delta)\,,$$
  where $\beta = r\alpha + (r-1)\cdot 2a'$. The issue is that the ``exit'' of $\iota_i$ does not align with the ``entry'' of $\iota_{i+1}$. More precisely, consider the sets $\iota_i((\alpha-\eta,\alpha)\times(0,\eta)^{2n-1})$ and $\iota_{i+1}((0,\eta)\times(0,\eta)^{2n-1})$. Their projections onto $0\times\R^{2n-1}\subset \R^{2n}$ are $(2ia'-\eta,2ia')^{2n-1}$ and $(2ia',2ia'+\eta)^{2n-1}$, respectively, so they are shifted relative to one another by the vector $(\eta,\dots,\eta)$, and thus the images of $\iota_i,\iota_{i+1}$ cannot be simply connected by a box to produce the embedding of a longer polydisk.

  To remedy this, we will shift each $B_{w^{(i)}}(a')$ for $i=1,\dots,r$ by the vector
  $$s_i:=-(0,(i-1)\eta,\dots,(i-1)\eta)\,,$$
  adjust the $\iota_i$ accordingly, and then the resulting embeddings will align, and we will be able to connect them by boxes. In detail, let $B_i':=T_{s_i}\big(B_{w^{(i)}}(a')\big)$ for $i=1,\dots,r$. Define
  $$\iota_i':=T_{s_i}\circ\iota_i \fc (0,\alpha)\times(0,\eta)^{2n-1}\hookrightarrow B_i'\quad\text{for }i=1,\dots,r\,.$$
  Note that now for $i=1,\dots,r-1$ the ``exit'' of $\iota_i'$ is aligned with the ``entry'' of $\iota_{i+1}'$. More precisely, the projections onto $0\times\R^{2n-1}\subset \R^{2n}$ of the sets $\iota_i'((\alpha-\eta,\alpha)\times(0,\eta)^{2n-1})$ and $\iota_{i+1}'((0,\eta)\times(0,\eta)^{2n-1})$ coincide. It follows that we can define a symplectic embedding
  $$\wh\iota \fc (0,\beta)\times(0,\eta)^{2n-1} \hookrightarrow\R^{2n}\quad\text{as follows:}$$
  $$\wh\iota(z)=\iota_i'(z)\quad \text{if}\quad u \in ((i-1)(\alpha+2a'),(i-1)(\alpha+2a')+\alpha)\;\text{for }i\in\{1,\dots,r\}
  \,,\text{ and}$$
  \begin{multline*}
    \wh\iota(z) = z+i(2a'-\alpha,2a'-\eta,\dots,2a'-\eta) \quad \text{if}\\ u \in (i\alpha+(i-1)\cdot2a'-\eta,i\alpha+i\cdot2a'+\eta)\text{ for }i\in\{1,\dots,r-1\}\,.
  \end{multline*}
  The image of $\wh\iota$ is schematically depicted in Figure \ref{fig:tubes_cubes}. Our condition \eqref{eqn:estimate_eta} on $\eta$ ensures that $\im\wh\iota \subset \bigcup_{i=1}^r B_{w^{(i)}}(a) \cup \bigcup_{i=1}^{r-1}N_i(\delta)$, since $\eta < \frac{\delta}{r-1}$. We also make the crucial observation that, thanks to Theorem \ref{thm:polydisk_into_cube}, item (i), for each $i\in\{1,\dots,r\}$, $\im\wh\iota$ contains the set
  $$B_{w^{(i)}}\big(a'-\max(2,r-1)\eta\big)\setminus\Sigma_\eta\,,$$
  where $\Sigma_\eta = \Sigma\big((\tfrac\eta2,\dots\tfrac\eta2),\eta\big)$ is the union of integer coordinate planes, scaled by a factor of $\eta$. By our choice of $\eta$, $a'-\max(2,r-1)\eta > a''$, meaning $B_{w^{(i)}}\big(a'-\max(2,r-1)\eta\big) \supset \ol B_{w^{(i)}}(a'')$.

  \begin{figure}
  \centering
  \includegraphics[width=12cm]{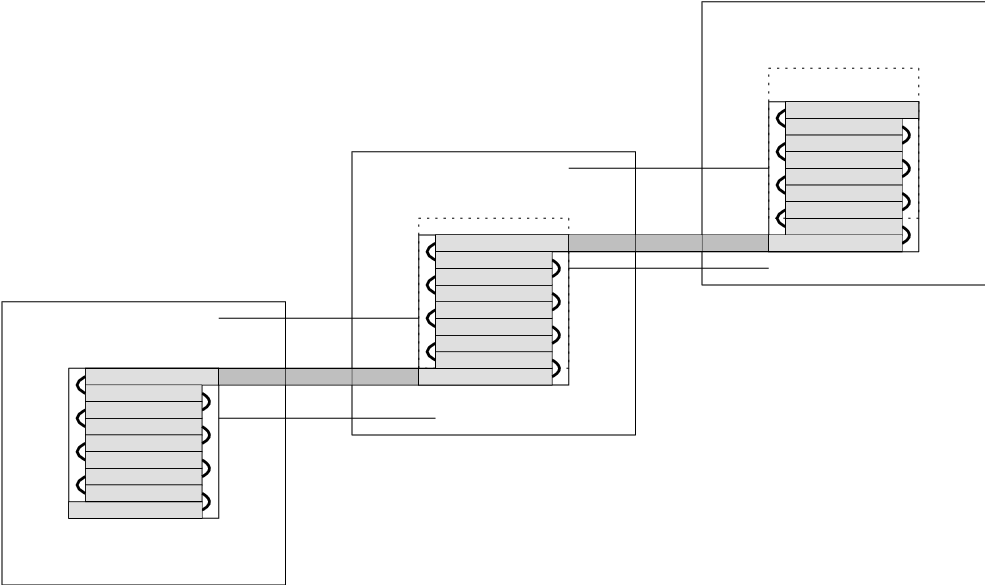}
  \caption{An illustration of the image of $\wh\iota$; the light shaded areas are the images of the corrected embeddings $\iota'_i$; the dark shaded areas are the connecting boxes. The squares with dotted boundary are the original cubes $B_{w^{(i)}}(a')$}\label{fig:tubes_cubes}
\end{figure}

  \tb{Step 6.} In this step we adjust the embedding we obtained in the previous step so that its image captures most of the measure $\mu$. Recall the embedding $\iota''$ from Step 3. In Step 5 we constructed an embedding $\wh\iota$ of a polydisk into the domain of $\iota''$, the idea being to eventually compose the two. This composition indeed yields an embedding of a polydisk into $M$, however, in addition, we need to control the measure of its image, and in this step we fix $\wh\iota$ to achieve this.

  Let
  $$\nu:=((\iota'')^{-1})_*(\mu|_{\im\iota''})$$
  be the measure induced on $\R^{2n}$ by pulling back $\mu$ via $\iota''$, where we view $(\iota'')^{-1}$ as a map $\im\iota''\to\R^{2n}$. By construction, the image of $\wh \iota$ contains the set
  $$\bigcup_{i=1}^r\ol B_{w^{(i)}}(a'')\setminus \Sigma_\eta\,.$$
  By our assumptions, $\nu\big(\bigcup_{i=1}^r\ol B_{w^{(i)}}(a'')\big)>1-\epsilon$. However, \emph{a priori}, $\nu(\Sigma_\eta)$ is unknown, and thus we cannot guarantee that $\nu\big(\bigcup_{i=1}^r\ol B_{w^{(i)}}(a'')\setminus \Sigma_\eta\big) > 1-\epsilon$, as claimed in the statement of the theorem. Let $\xi \in \R^{2n}$ be such that
  $$\nu\big(T_\xi(\Sigma_\eta)\big)=0\,,$$
  and such that $\|\xi\|<\min\big(a'-a''-(r-1)\eta,a'-a''-2\eta\big)$. Such $\xi$ exists thanks to Lemma \ref{lem:Sigma_y_a_at_most_countable}. It follows that
  $$\textstyle\nu\Big(\bigcup_{i=1}^r\ol B_{w^{(i)}}(a'')\setminus T_\xi(\Sigma_\eta)\Big) = \nu\Big(\bigcup_{i=1}^r\ol B_{w^{(i)}}(a'')\Big) > 1-\epsilon\,.$$
  Let us denote
  $$B_i'' = B_{w^{(i)}}(a'-\max(2,r-1)\eta)\,.$$
  At the end of Step 5 we noted that $B_i'' \supset \ol B_{w^{(i)}}(a'')$ for each $i$, and that
  $$\im\wh\iota\supset \bigcup_{i=1}^rB_i''\setminus\Sigma_\eta\,.$$
  By our choice of $\xi$, we also have $T_\xi(B_i'')\supset \ol B_{w^{(i)}}(a'')$, whence
  \begin{equation}\label{eq:cube_containment_correcting_wh_iota}
    \bigcup_{i=1}^rT_\xi(B_i'')\setminus T_\xi(\Sigma_\eta) \supset \bigcup_{i=1}^r \ol B_{w^{(i)}}(a'')\setminus T_\xi(\Sigma_\eta)\,.
  \end{equation}
  Again, by the choice of $\xi$, for $i=1,\dots,r$ there exist closed cubes $C_i,C_i'$ such that $C_i\supset B_i''\cup T_\xi(B_i'')$, such that $C_i$ is contained in the interior of $C_i'$, and such that $C_i'\subset B_{w^{(i)}}(a)$. Let $K$ be a Hamiltonian on $\R^{2n}$ which is obtained by multiplying the function $-\langle \sqrt{-1}\xi,\cdot\rangle$ (here $\sqrt{-1}$ is the imaginary unit) by a cutoff function which on each $C_i$ is identically $1$, and which vanishes outside $\bigcup_{i=1}^r C_i'$. It then follows that the restriction of $\phi_K$ to each $B_i''$ coincides with the translation $T_\xi$. Note as well that $\phi_K$ is a symplectomorphism with compact support inside $\bigcup_{i=1}^r B_{w^{(i)}}(a) \cup \bigcup_{i=1}^{r-1}N_i(\delta)$. Let $\phi$ be the restriction of $\phi_K$ to this set. The corrected embedding now is
  $$\phi \circ \wh\iota \fc (0,\beta)\times(0,\eta)^{2n-1} \hookrightarrow \bigcup_{i=1}^r B_{w^{(i)}}(a) \cup \bigcup_{i=1}^{r-1}N_i(\delta)\,.$$
  By construction, its image contains $\bigcup_{i=1}^r T_\xi(B_i'')\setminus T_\xi(\Sigma_\eta)$, and thus by equation \eqref{eq:cube_containment_correcting_wh_iota},
  $$\im\phi\circ\wh\iota \supset \bigcup_{i=1}^r\ol B_{w^{(i)}}(a'')\setminus T_\xi(\Sigma_\eta)\,,\quad\text{whence}$$
  $$\textstyle\nu\big(\im(\phi\circ\wh\iota)\big) \geq \nu\Big(\bigcup_{i=1}^r\ol B_{w^{(i)}}(a'')\setminus T_\xi(\Sigma_\eta)\Big) = \nu\Big(\bigcup_{i=1}^r\ol B_{w^{(i)}}(a'')\Big) > 1-\epsilon\,.$$

\begin{figure}
  \centering
  \includegraphics[width=8cm]{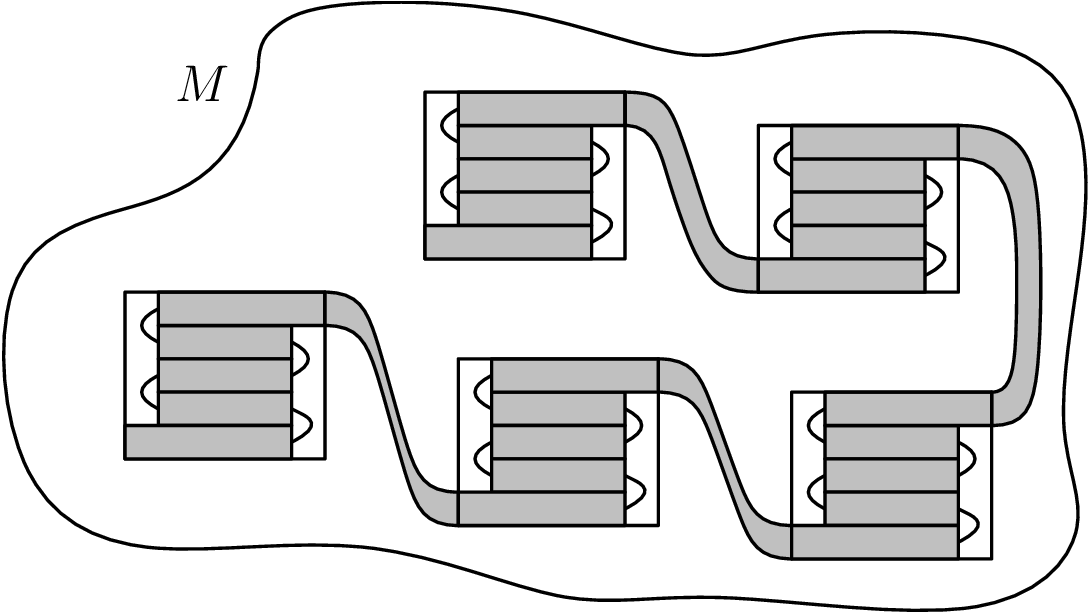}
  \caption{The image of the embedding $\iota''$ inside the manifold $M$}\label{fig:embedded_polydisk_M}
\end{figure}

  \tb{Step 7.} Finally we assemble the above embeddings to obtain the following:
  $$\iota '' \circ \phi\circ\wh\iota \fc (0,\beta)\times (0,\eta)^{2n} \hookrightarrow M\,.$$
  The image of the resulting embedding is schematically depicted in Figure \ref{fig:embedded_polydisk_M}. From the previous step it follows that
  $$\mu\big(\im(\iota''\circ\phi\circ\wh\iota)\big) = \nu\big(\im(\phi\circ\wh\iota)\big)>1-\epsilon\,.$$
  The proof of the theorem is complete.
\end{proof}

\section{Proof of Main result}\label{s:pf_main_result}

Here we prove the Main result, Theorem \ref{thm:main_result}, which states that given a closed connected symplectic manifold $(M,\omega)$ of dimension $\geq 4$ and with $H^1(M;\Z)=0$, the only symplectic Aarnes quasi-states on $M$ are delta-measures. As we have outlined in the introduction, the strategy of the proof is as follows. Fix $\mu \in \cP_0(M)$ (recall Notation \ref{not:P_0_X}), and let $\zeta$ and $\tau$ be the corresponding Aarnes quasi-state and the topological measure representing it, respectively. The first step in the proof is Theorem \ref{thm:main_embedding}, which yields a symplectic embedding $\iota \fc P \hookrightarrow M$, where $P$ is a polydisk and $\mu(\iota(P))>\frac 12$. Theorem \ref{thm:main_embedding} yields such an embedding whose image has $\mu$-measure arbitrarily close to $1$, but all we need here is for it to be $>\frac 12$.

In the second step, performed in Section \ref{ss:involutive_map_fibers}, we construct an involutive map $\Phi \fc M \to \R^n$, where $n=\frac 12\dim M$. It has the property that there is a slightly smaller polydisk $P'\subset P$, such that we still have $\mu(\iota(P'))>\frac 12$, and such that the fiber components of $\Phi$ are $M\setminus \iota(P')$ and the isotropic tori $\iota(T(\alpha))$, whose union is $\iota(P')$, where the $T(\alpha)$ were defined in Example \ref{ex:Lag_tori_T_alpha}. Thus the special fiber component $C_\Phi$ of $\Phi$ with respect to $\tau$ (see Definition \ref{def:special_fiber_cpt}) is either $M\setminus \iota(P')$ or one of these tori. Each fiber component is solid, thanks to Example \ref{ex:solid_sets} (\emph{this is where we use the assumption that} $\dim M \geq 4$), whence $\mu(C_{\phi})\geq\frac12$, thanks to Lemma \ref{lem:solid_cpt_sh_meas_more_half}, which forces $C_\Phi$ to be one of the tori, since $\mu(M\setminus \iota(P')) < \frac 12$ by construction.

We are now in the situation where there is an isotropic torus $T'\subset M$ with $\mu(T')\geq \frac12$. This torus is either itself Lagrangian, or can easily be shown to be contained in a Lagrangian torus. In any case, we obtain a Lagrangian torus $T\subset M$ with $\mu(T)\geq \frac12$. In the third and final step, which appears in Section \ref{ss:sh_Lag_has_sh_point}, we use a combination of set- and measure-theoretic arguments to show that the symplecticity of $\zeta$ forces the existence of a point $q \in T$ with $\mu(\{q\})\geq \frac12$, which, thanks to Corollary \ref{cor:Aarnes_delta}, means that $\tau=\delta_q$.

\subsection{The involutive map and its fibers}\label{ss:involutive_map_fibers}

In this section we will use the map $\Phi_{\std} \fc \C^n \to \R^n$, $\Phi_{\std}(z) = (|z_1|^2,\dots,|z_n|^2)$, which is involutive thanks to Example \ref{ex:invol_maps}, item (iv). For the remainder of this section we fix a smooth even function $\chi \fc \R\to[0,1]$ which vanishes exactly on $\R\setminus(-1,1)$, and which is strictly decreasing on $[0,1]$, and therefore strictly increasing on $[-1,0]$. For $r > 0$ we will use the scaled version $\chi_r=\chi(\frac{\cdot}{r})$.

For $a=(a_1,\dots,a_n) \in (0,\infty)^n$ we define the following map $\rho_a \fc \R^n\to \R^n$:
\begin{equation}\label{eqn:def_rho_a}
  \rho_a(t) = \chi_{a_1}(t_1)\cdot\ldots\cdot\chi_{a_n}(t_n)(a_1^2-t_1^2,\dots,a_n^2-t_n^2)\,.
\end{equation}
Note that it is smooth and that it vanishes exactly outside $\prod_{j=1}^{n}(-a_j,a_j)$. Combining this with $\Phi_{\std}$, we obtain $\wt\Phi_a \fc \C^n \to \R^n$, $\wt\Phi_a = \rho_a\circ\Phi_{\std}$, which is again involutive thanks to Example \ref{ex:invol_maps}, item (vii).

Recall that we have denoted $D(r)=\{z \in \C\,|\,|z|^2<r\}$ for $r>0$ and $P(b) = \prod_{j=1}^{n}D(b_j)$ for $b=(b_1,\dots,b_n)\in(0,\infty)^n$. Another notation we will use is $E(b) = \prod_{j=1}^{n}[0,b_j)$. Note that $\Phi_{\std}(P(b)) = E(b)$.

Let $(M,\omega)$ be a closed connected symplectic manifold, and let $\iota\fc P(b) \hookrightarrow M$ be a symplectic embedding. Pick $a_j \in (0,b_j)$ for all $j$, put $a=(a_1,\dots,a_n)$, and define $\Phi_{\iota;a} \fc M \to \R^n$ by
\begin{equation}\label{eqn:def_Phi_iota_a}
  \Phi_{\iota;a}|_{M\setminus\im\iota}\equiv 0\quad\text{and}\quad \Phi_{\iota;a}|_{\im \iota}=\wt\Phi_a\circ\iota^{-1}\,.
\end{equation}

\begin{lemma}\label{lem:Phi_iota_a_smooth_invol}
  The map $\Phi_{\iota;a}$ is smooth and involutive.
\end{lemma}
\begin{proof}
  For the smoothness consider the open cover of $M$ by the sets $\im\iota$ and $M\setminus(\iota(\ol{P(a)}))$: on $\im\iota$, $\Phi_{\iota;a}$ is smooth as the composition of smooth maps, while on $M\setminus(\iota(\ol{P(a)}))$ it is the constant map with value $0$, in particular it is smooth. For the involutivity note that it is likewise a local property, therefore it suffices to check it on the same two open subsets. The zero map is obviously involutive. Lastly, $\wt\Phi_a\circ\iota^{-1}$ is involutive as the composition of the involutive map $\wt\Phi_a$ and the symplectomorphism $\iota^{-1} \fc \im\iota \to P(b)$, see Example \ref{ex:invol_maps}, item (vi).
\end{proof}

We will now calculate all the fiber components of $\Phi_{\iota;a}$. Recall the notation of Example \ref{ex:Lag_tori_T_alpha}: $T(\alpha) = \Phi_{\std}^{-1}(\alpha)$ for $\alpha \in[0,\infty)^n$. Then we have
\begin{prop}\label{prop:Phi_iota_fiber_finite}
Each fiber of\; $\Phi_{\iota;a}$ is connected, in particular $\Phi_{\iota;a}$ is fiber-finite. Each fiber (component) is either $M\setminus \iota(P(a))$, or $\iota(T(\alpha))$ for some $\alpha \in E(a)$.
\end{prop}
\noindent We will first prove the following lemma detailing the relevant properties of $\rho_a$.
\begin{lemma}\label{lem:rho_a_finite_fibers}
  For $c \in \R^n$, we have the following:
  $$\rho_a^{-1}(c) = \begin{cases}\varnothing\,,&c\notin\im\rho_a\,, \\ \R^n\setminus \prod_{j=1}^{n}(-a_j,a_j)\,,& c=0 \\ \prod_{j=1}^{n}\{-\alpha_j,\alpha_j\}\,,&\text{otherwise}\,,\end{cases}$$
  where for each $j$, $\alpha_j \in [0,a_j)$ is a number depending on $c$.
\end{lemma}
\noindent Note that $\{-\alpha_j,\alpha_j\}=\{0\}$ if $\alpha_j=0$.
\begin{proof}
  Clearly $c \notin \im\rho_a$ if and only if $\rho_a^{-1}(c) = \varnothing$. Thus two cases remain, namely (i) $c=0$ and (ii) $c \in \im \rho_a\setminus\{0\}$.

  \tb{Case (i):} $c=0$. As we mentioned after the definition of $\rho_a$, equation \eqref{eqn:def_rho_a}, it vanishes exactly outside $\prod_{j=1}^{n}(-a_j,a_j)$.

  \tb{Case (ii):} We have $c\in\im\rho_a\setminus\{0\}$. Any solution of $\rho_a(t) = c$ satisfies $t \in \prod_{j=1}^{n}(-a_j,a_j)$, and the equation is equivalent to the following system, where $c=(c_1,\dots,c_n)$:
  \begin{equation}\label{eqn:system_rho_equals_c}
    \forall j: \quad \chi_{a_1}(t_1)\dots\chi_{a_n}(t_n)(a_j^2-t_j^2) = c_j\,.
  \end{equation}
  Since for each $j$ we have $|t_j|<a_j$, from our choice the function $\chi$ it follows that $\chi_{a_j}(t_j) > 0$, and thus the left-hand side of the equation is positive, which implies that $c_j>0$ for all $j$. In case $n=1$ we have the sole equation $\chi_{a_1}(t_1)(a_1^2-t_1^2)=c_1$. The function on the left-hand side is strictly decreasing on $[0,a_1]$, and strictly increasing on $[-a_1,0]$. It follows that this equation has either the unique solution $t_1=0$ in case $c_1=\chi_{a_1}(0)a_1^2$, or exactly two solutions, which moreover have the same absolute value.

  In case $n>1$ we can divide all the equations, starting from the second one, by the first one, which results in
  $$\frac{a_j^2-t_j^2}{a_1^2-t_1^2}=\frac{c_j}{c_1}\,.$$
  Denoting $d_j=\frac{c_j}{c_1}$, we can solve this equation for $t_j$ to get
  \begin{equation}\label{eqn:t_j_in_terms_t_1}
    t_j=\pm\sqrt{a_j^2-d_j(a_1^2-t_1^2)}\,.
  \end{equation}
  Substituting this into the first equation of the system \eqref{eqn:system_rho_equals_c}, and remembering that the $\chi_{a_j}$ are even, we obtain the following equation for $t_1$:
  \begin{equation}\label{eqn:equation_t_1}
    \chi_{a_1}(t_1)\chi_{a_2}\Big(\sqrt{a_2^2-d_2(a_1^2-t_1^2)}\Big)\dots\chi_{a_n}\Big(\sqrt{a_n^2-d_j(a_1^2-t_1^2)}\Big)(a_1^2-t_1^2)=c_1\,.
  \end{equation}
  Consider the expression $\sqrt{a_j^2-d_j(a_1^2-t_1^2)}$ as a function of $t_1\geq 0$. It is well-defined as long as $t_1\geq\sqrt{a_1^2 - a_j^2/d_j}$, and there it is a strictly increasing function of $t_1$. It follows that, as long as we are considering $t_1\geq0$, the left-hand side of \eqref{eqn:equation_t_1} is defined on some interval of the form $t_1\geq\const$, and it is a strictly decreasing function of $t_1$ there until some larger value of $t_1$ starting from which it vanishes. Therefore the above equation has exactly one solution for $t_1\geq 0$. Since the function is even, the only other solution is obtained from the one we have just obtained by negating it. It follows that the above equation always has either exactly one solution, namely $t_1=0$ or exactly two solutions, which are equal in absolute value. Finally, since $t_j$ can be recovered from $t_1$ by \eqref{eqn:t_j_in_terms_t_1}, we see that the solution set of the system has the form asserted in the formulation of the lemma. The proof is complete.
\end{proof}

\begin{proof}[Proof of Proposition \ref{prop:Phi_iota_fiber_finite}.]
  We abbreviate $\Phi=\Phi_{\iota;a}$. We have
  $$\Phi^{-1}(0) = (\Phi|_{M\setminus\im\iota})^{-1}(0) \cup (\Phi|_{\im\iota})^{-1}(0)\,.$$
  Since $\Phi|_{M\setminus\im\iota}\equiv 0$, we have $(\Phi|_{M\setminus\im\iota})^{-1}(0) = M\setminus\im\iota$. Also $\Phi|_{\im\iota} = \rho_a\circ\Phi_{\std}\circ\iota^{-1}$. To calculate $(\Phi|_{\im\iota})^{-1}(0)$, let us use the more precise formula
  $$\Phi|_{\im\iota} = \rho_a|_{E(b)}\circ\Phi_{\std}|_{P(b)}\circ\iota^{-1}\,,$$
  whence
  $$(\Phi|_{\im\iota})^{-1}(0) = (\rho_a|_{E(b)}\circ\Phi_{\std}|_{P(b)}\circ\iota^{-1})^{-1}(0) = \iota((\Phi_{\std}|_{P(b)})^{-1} ((\rho_a|_{E(b)})^{-1}(0)))\,.$$
  We have, by Lemma \ref{lem:rho_a_finite_fibers}:
  $$\textstyle(\rho_a|_{E(b)})^{-1}(0) = {E(b)}\cap\big(\R^n\setminus\prod_{j=1}^{n}(-a_j,a_j)\big) = E(b)\setminus E(a)\,,$$
  which implies
  $$(\Phi_{\std}|_{P(b)})^{-1}((\rho_a|_{E(b)})^{-1}(0)) = (\Phi_{\std}|_{P(b)})^{-1}(E(b)\setminus E(a)) =  P(b)\setminus P(a)\,,$$
  that is
  $$\iota((\Phi_{\std}|_{P(b)})^{-1}((\rho_a|_{E(b)})^{-1}(0))) = \iota(P(b)\setminus P(a)) = \im\iota\setminus\iota(P(a))\,.$$
  In total we obtain
  $$\Phi^{-1}(0) = (M\setminus\im\iota)\cup (\im\iota\setminus\iota(P(a))) = M\setminus\iota(P(a))\,.$$
  Since the set on the right is connected thanks to Example \ref{ex:solid_sets}, it is the unique component of $\Phi^{-1}(0)$.

  Now assume $c \in \im\Phi\setminus\{0\}$. It follows that $\Phi^{-1}(c) = (\Phi|_{\im\iota})^{-1}(c)$, thus
  $$\Phi^{-1}(c) = \iota((\Phi_{\std}|_{P(b)})^{-1}(\rho_a|_{E(b)})^{-1}(c))\,.$$
  Since $c \neq 0$, by Lemma \ref{lem:rho_a_finite_fibers}, $\rho_a^{-1}(c) = \prod_{j=1}^{n}\{-\alpha_j,\alpha_j\}$ for some numbers $\alpha_j \in [0,a_j)$. In particular $(\rho_a|_{E(b)})^{-1}(c) = (\alpha_1,\dots,\alpha_n)$. Denoting $\alpha=(\alpha_1,\dots,\alpha_n)$, we obtain
  $$\Phi^{-1}(c) = \iota (\Phi_{\std}^{-1}((\rho_a|_{E(b)})^{-1}(c))) = \iota(\Phi_{\std}^{-1}(\alpha)) = \iota(T(\alpha))\,.$$
  The proof of the proposition is complete.
\end{proof}

\begin{prop}\label{prop:sympl_Aarnes_QS_sh_Lag_torus}
  Let $(M,\omega)$ be a closed connected symplectic manifold of dimension $\geq 4$ with $H^1(M;\Z)=0$ and let $\mu \in \cP_0(M)$. Let $\zeta$ be the corresponding Aarnes quasi-state and let $\tau$ be the topological measure representing it. If $\zeta$ is symplectic, then there exists a Lagrangian torus $T\subset M$ with $\mu(T)\geq \frac12$. In particular, $T$ is $\zeta$-superheavy.
\end{prop}
\begin{proof}
  By Theorem \ref{thm:main_embedding}, there exists a symplectic embedding $\iota \fc P(b) \hookrightarrow M$ for some $b\in(0,\infty)^n$, such that $\mu(\im\iota)>\frac 1 2$. Let $a\in E(b)$ be such that $\mu(\iota(P(a)))>\frac 1 2$, which exists by Lemma \ref{cor:approx_open_set_cpts}. Let $\Phi:=\Phi_{\iota;a} \fc M \to \R^n$ be the map constructed above, see equation \eqref{eqn:def_Phi_iota_a}. By Proposition \ref{prop:Phi_iota_fiber_finite}, $\Phi$ is fiber-finite, therefore we have the corresponding special fiber component $C_\Phi$, which satisfies $\tau(C_\Phi)=1$, see Lemma \ref{lem:fiber_finite_invol_maps_special_fiber_cpt}. Again by Proposition \ref{prop:Phi_iota_fiber_finite}, the fiber components of $\Phi$ are the set $M\setminus\iota(P(a))$ and the isotropic tori $\iota(T(\alpha))$ with $\alpha \in E(a)$.

  Any such connected component is solid, thanks to Example \ref{ex:solid_sets}. \emph{Note that it is here that we use the assumption on the dimension of $M$}. Indeed, as Example \ref{ex:solid_sets} shows, an isotropic torus in a connected symplectic manifold of dimension $\geq4$ is solid. This is clearly false in dimension two, since there are disconnecting circles on any surface.

  Lemma \ref{lem:solid_cpt_sh_meas_more_half} then implies that $\mu(C_\Phi)\geq\frac 1 2$, and in particular $C_\Phi$ cannot be $M\setminus\iota(P(a))$. It follows that $C_\Phi = \iota(T(\alpha))$ for some $\alpha \in E(a)$.

  If for all $j$ we have $\alpha_j>0$, put $T=\iota(T(\alpha))$, and note that it is Lagrangian. Otherwise $T(\alpha) = \prod_{j=1}^{n} \{z \in \C\,|\,|z|^2=\alpha_j\}$. If $\alpha_j>0$, put $T_j=\{z \in \C\,|\,|z|^2=\alpha_j\}$, otherwise define
  $$\textstyle T_j=\big\{z\in\C\;|\;\big|z-\frac{\sqrt{a_j}}{4}\big|=\frac{\sqrt{a_j}}{4}\big\}$$
  and note that $0\in T_j\subset D(a_j)$. Therefore $T(\alpha)\subset \prod_{j=1}^{n}T_j\subset P(a)$ and we can define $T=\iota\big(\prod_{j=1}^{n}T_j\big)$, which is a Lagrangian torus.

  In any case we have constructed a Lagrangian torus $T$ with $C_\Phi \subset T$, which implies that $\mu(T)\geq\frac12$, as required. Lemma \ref{lem:solid_cpt_sh_meas_more_half} then yields $\tau(T) = 1$, whence $T$ is $\zeta$-superheavy thanks to Lemma \ref{lem:char_sh_sets_tau_is_one}.
\end{proof}

\subsection{Superheavy Lagrangian torus contains superheavy point}\label{ss:sh_Lag_has_sh_point}

This section is dedicated to the proof of Theorem \ref{thm:sh_Lag_torus_contains_sh_point}, which, together with Proposition \ref{prop:sympl_Aarnes_QS_sh_Lag_torus} and Corollary \ref{cor:Aarnes_delta}, implies the main result, Theorem \ref{thm:main_result}. Recall that Theorem \ref{thm:sh_Lag_torus_contains_sh_point} states:

\emph{Let $(M,\omega)$ be a closed connected symplectic manifold such that $\dim M \geq 4$ and $H^1(M;\Z)=0$, and let $\zeta$ be the Aarnes quasi-state on $M$ corresponding to a measure $\mu$ as in Theorem \ref{thm:Aarnes_constr}, and assume that $\zeta$ is symplectic; if $T\subset M$ is a $\zeta$-superheavy Lagrangian torus, then there is $z_0 \in T$ with $\mu(\{z_0\}) \geq \frac 12$.}

\begin{rem}\label{rem:superheavy_torus_meas_at_leaste_half}
  Note that since $\dim M \geq 4$, a Lagrangian torus is necessarily solid, which means, thanks to Lemmas \ref{lem:char_sh_sets_tau_is_one} and \ref{lem:solid_cpt_sh_meas_more_half}, that $T$, being $\zeta$-superheavy, satisfies $\mu(T) \geq \frac12$.
\end{rem}

The proof appears at the end of the section. For it we will need preparatory results. Let us summarize their main ideas. All the ``action'' happens in a neighborhood of the given superheavy Lagrangian torus $T$. Since $T$ is a torus, it has a neighborhood which is symplectomorphic to $T\times B$, where $B$ is an open neighborhood of zero in $\R^n\cong T^*_qT$, see Remark \ref{rem:Weinstein_nbd_thm}. Using this splitting, we will define an involutive map $\Phi$ on $M$ whose fiber components are the tori $T\times\{p\}$, where $p$ ranges over $B$, and the complement $M\setminus(T\times B)$. Since $T$ is superheavy by assumption, it is the special fiber component of $\Phi$. Now we use the flexibility afforded by the multitude of symmetries in symplectic geometry given by Hamiltonian diffeomorphisms. In particular, for each $f \in C^\infty(T)$, whose $C^1$-norm is small enough, we will produce another involutive map $\Phi_f$ on $M$ whose special fiber component $T_f$ is the graph of a closed $1$-form on $T$ of the form $\kappa(f)+df$ for some constant $\kappa(f) \in \R^n$. Remark \ref{rem:superheavy_torus_meas_at_leaste_half} applies to $T_f$ as well, and yields $\mu(T_f)\geq \frac 12$.

We then use the freedom of choice of $f$ to deduce, based on set- and measure-theoretic arguments, that $T$ contains a point of $\mu$-measure $\geq \frac12$, completing the proof of Theorem \ref{thm:sh_Lag_torus_contains_sh_point}, and therefore of the main result Theorem \ref{thm:main_result}.

Let us now describe all of this in detail. We identify $\T^n=\R^n/\Z^n$ and $T^*\T^n = \T^n(q)\times\R^n(p)$ which is equipped with the canonical symplectic form $dp\wedge dq$. The projection onto the second factor $\pr \fc T^*\T^n=\T^n\times\R^n\to\R^n$ is involutive, see Example \ref{ex:invol_maps}, item (v).

\paragraph*{Notation.} In this section we denote the action of $\R^n$ by translation on the second factor of $\T^n\times\R^n$ by $\xi+(q,p)$, that is $\xi+(q,p)=(q,p+\xi)$. Moreover, for a $1$-form $\sigma$ on $\T^n$ we let $\Gamma_\sigma\subset T^*\T^n$ be its graph: $\Gamma_\sigma = \{(q,\sigma_q)\,|\,q\in \T^n\}$. We will also use the uniform norm $\|\sigma\|$ of $\sigma$ relative to the $\ell^\infty$-metric on the fiber $\R^n$, that is $\|\sigma\|=\sup\{\|\sigma_q\|\,|\,q\in \T^n\}$.

We fix a diffeomorphism $T \cong \T^n=\R^n/\Z^n$. By the Weinstein neighborhood theorem, there exist $\epsilon'>0$ and a neighborhood of $T$ which is symplectomorphic to $T \times B_0(\epsilon') = \T^n(q)\times (-\epsilon',\epsilon')^n(p) \subset T^*\T^n$. We will identify such a neighborhood with $T \times B_0(\epsilon')$ and $T$ with the zero section $\T^n\times \{0\}=\T^n$. Fix $\epsilon \in (0,\epsilon')$, let $a=(\epsilon,\dots,\epsilon) \in (0,\infty)^n$ and define $\Phi \fc M \to \R^n$ by
$$\Phi(z)|_{M\setminus (T\times B_0(\epsilon))}\equiv 0\,,\quad \Phi|_{T\times B_0(\epsilon)} = \rho_a\circ \pr\,,$$
where $\rho_a$ is defined by \eqref{eqn:def_rho_a}. The necessary properties of $\Phi$ are summarized as follows.
\begin{lemma}\label{lem:invol_map_Phi_torus_fiber_finite}
  The map $\Phi$ is smooth, involutive, and fiber-finite. Moreover, its fiber components are the set $M\setminus(T\times B_0(\epsilon))$ and the Lagrangian tori $T\times \{p\}$ for $p \in B_0(\epsilon)$.
\end{lemma}
\begin{proof}
  The smoothness and involutivity are checked in the same way as in the proof of Lemma \ref{lem:Phi_iota_a_smooth_invol}. Regarding the fiber components of $\Phi$, note that $\rho_a$ vanishes exactly outside $B_0(\epsilon)$, therefore
  $$\Phi^{-1}(0) = (M\setminus(T\times B_0(\epsilon)))\cup \pr^{-1}(\underbrace{(\rho_a|_{B_0(\epsilon)})^{-1}(0)}_{=\varnothing}) = M\setminus(T\times B_0(\epsilon))\,.$$
  That this set is connected follows from Example \ref{ex:solid_sets}, item (iii). For $c \in \im\Phi\setminus\{0\}$, $\rho_a^{-1}(c)$ is a finite set by Lemma \ref{lem:rho_a_finite_fibers}, and thus
  $$\Phi^{-1}(c) = \pr^{-1}(\rho_a^{-1}(c))=\bigsqcup_{p\in\rho_a^{-1}(c)}\pr^{-1}(p) = \bigsqcup_{p\in\rho_a^{-1}(c)}T\times\{p\}\,,$$
  as claimed.
\end{proof}
Now for any symplectomorphism $\psi$ of $M$, the map $\Phi\circ\psi$ is still involutive, thanks to Example \ref{ex:invol_maps}, item (vi). Lemma \ref{lem:fiber_finite_invol_maps_special_fiber_cpt} then yields the special fiber component $C_{\Phi\circ\psi}$. Since $\psi$ is a homeomorphism, the fiber components of $\Phi\circ\psi$ are the set $\psi^{-1}(M\setminus(T\times B_0(\epsilon))) = M\setminus\psi^{-1}(T\times B_0(\epsilon))$ and the tori $\psi^{-1}(T\times\{p\})$ for $p\in B_0(\epsilon)$. Since $\psi$ is a symplectomorphism, these tori are likewise Lagrangian.

Now assume that $\psi(T\times B_0(\epsilon)) = T\times B_0(\epsilon)$. It follows that the fiber components of $\Phi\circ\psi$ are $M\setminus(T\times B_0(\epsilon))$ and the tori $\psi^{-1}(T\times\{p\})$ for $p \in B_0(\epsilon)$. By Remark \ref{rem:special_fiber_cpts_intersect}, the special fiber components $C_{\Phi\circ\psi}$ and $C_\Phi=T$ must intersect. It follows that $C_{\Phi\circ\psi}$ must be of the form $\psi^{-1}(T\times \{p\})$, and, of course, it must intersect $T$.

Let $\theta \in C^\infty(M,[0,1])$ be a smooth function with $\theta|_{T\times \ol B_0(2\epsilon/3)}\equiv 1$ and $\theta|_{M\setminus (T\times B_0(\epsilon))}\equiv 0$. For $f \in C^\infty(T)$ define $H_f=\theta\cdot\pi^*f$ on $T\times B_0(\epsilon)$, where $\pi \fc T^*\T^n\to\T^n$ is the base projection, and extend it by zero to $M$. Let $\phi_f^t:=\phi_{H_f}^t$ and $\phi_f:=\phi_{H_f}$. We record the following fundamental property of $\phi^t_f$:
\begin{lemma}\label{lem:flow_phi_f_same_as_pullback_f}
  Assume that $f \in C^\infty(T)$ satisfies $\|df\| \leq \frac{\epsilon}3$. Then for all $t \in [0,1]$ we have $\phi_f^t|_{T\times\ol B_0(\frac{\epsilon}3)}=\phi^t_{\pi^*f}=\big((q,p)\mapsto(q,p-td_qf)\big)$.
\end{lemma}
\begin{proof}
  This is a consequence of the following basic fact: if $X,Y$ are complete vector fields on a manifold $Q$, $K,L\subset Q$ are compact subsets such that $\bigcup_{t\in[0,1]}\phi_X^t(K)\subset L$, and $Y|_L\equiv X|_L$, then $\phi_{Y}^t|_K=\phi_X^t|_K$ for $t\in[0,1]$. The lemma follows by substituting $Q=T^*\T^n=T^*T$, $K=T\times\ol B_0(\epsilon/3)$, $L=T\times\ol B_0(2\epsilon/3)$, $X=X_{\pi^*f}$, $Y=X_{H_f}$, where $H_f$ is extended by zero to the whole of $T^*\T^n$.
\end{proof}

\begin{figure}
  \centering
  \includegraphics[width=12cm]{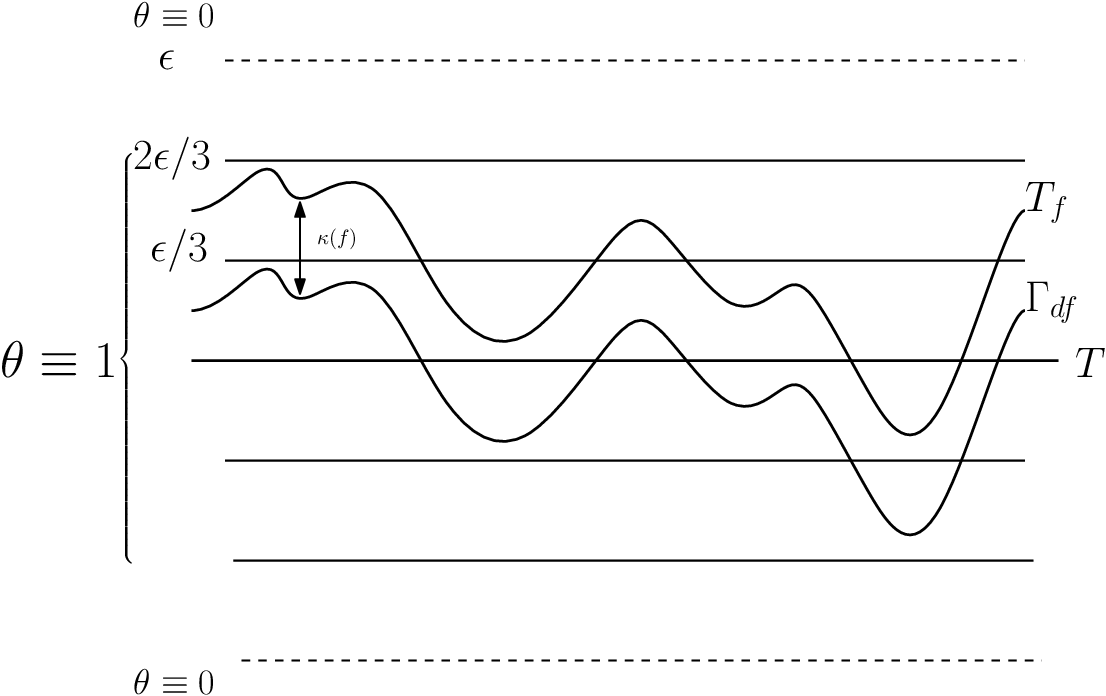}
  \caption{The torus $T$, the cutoff function $\theta$, the graph of $df$ and the special fiber $T_f$}\label{fig:lemma_4_8}
\end{figure}

\begin{lemma}\label{lem:T_f_is_Gamma_f_plus_kappa}
  For $f\in C^\infty(T)$ with $\|df\|\leq \frac{\epsilon}{3}$, the special fiber component $C_{\Phi\circ\phi_f}$ has the form $\kappa(f)+\Gamma_{df}$, where $\kappa(f) \in \ol B_0(\|df\|)$.
\end{lemma}
\noindent See Figure \ref{fig:lemma_4_8} for an illustration.

\begin{proof}
  Since $H_f$ vanishes outside $T\times B_0(\epsilon)$, $\phi_f$ preserves this set, therefore by the above discussion we have $$C_{\Phi\circ\phi_f}=\phi_f^{-1}(T\times\{p\})\qquad \text{for some }p \in B_0(\epsilon)\,.$$
  It must intersect $T$: $T\cap \phi_f^{-1}(T\times\{p\}) \neq \varnothing$, or equivalently $T\times\{p\}\cap \phi_f(T)\neq \varnothing$. However, thanks to Lemma \ref{lem:flow_phi_f_same_as_pullback_f}, we have $\phi_f(T) = \Gamma_{-df}$, therefore we conclude that $p \in \ol B_0(\|df\|)$. By the same lemma, $\phi_f^{-1}(T\times\{p\}) = \phi_{\pi^*f}^{-1}(T\times\{p\}) = p+\Gamma_{df}$, which completes the proof with $\kappa(f) = p$.
\end{proof}

\begin{notation}\label{not:special_fiber_cpt_f}
  For $f \in C^\infty(T)$ such that $\|df\| \leq \frac{\epsilon}{3}$ we let $T_f=\kappa(f)+\Gamma_{df}\subset T\times B_0(\epsilon)$ stand for the special fiber component $C_{\Phi\circ\phi_f}$.
\end{notation}

\begin{rem}\label{rem:T_f_contained_ball_2_norm_df}
   Note that $T_f\subset T\times \ol B_0(2\|df\|)\subset T\times \ol B_0(2\epsilon/3)$, since $\|\kappa(f)+df\|\leq 2\|df\|$.
\end{rem}

Given this preparation, we can now prove Theorem \ref{thm:sh_Lag_torus_contains_sh_point}, for which we need two additional lemmas.

\begin{lemma}\label{lem:point_half_mass}
  Let $(X,\mu)$ be a measure space of finite total measure, where singletons are measurable. If for each $\sigma \in (0,\frac12)$ there exists $x\in X$ with $\mu(\{x\})>\frac12-\sigma$, then there exists $x \in X$ with $\mu(\{x\}) \geq \frac 12$.
\end{lemma}
\begin{proof}
  For each $\sigma \in (0,\frac 12)$ let $X_\sigma = \{x \in X\,|\,\mu(\{x\})>\frac12-\sigma\}$. Then $(X_\sigma)_{\sigma\in(0,\frac 12)}$ is a collection of nonempty finite sets, which decreases as $\sigma$ decreases. Let $X_0 = \bigcap_{\sigma\in(0,\frac 12)}X_\sigma$. Since $X_0$ is a decreasing intersection of nonempty finite sets, it is itself nonempty. We claim that any $x \in X_0$ satisfies $\mu(\{x\}) \geq \frac12$. Indeed, since any $x \in X_0$ satisfies $x \in X_\sigma$ for all $\sigma\in(0,\frac 12)$, we have $\mu(\{x\}) > \frac 12-\sigma$. Taking the limit as $\sigma \to 0$, we arrive at the desired conclusion.
\end{proof}

In the next lemma, as well as in the proof of Theorem \ref{thm:sh_Lag_torus_contains_sh_point} immediately following it, we use additive notation for the group structure on $S^1=\R/\Z$.
\begin{lemma}\label{lem:set_reflection}
  Let $Q\subset S^1=\R/\Z$ be a subset which is at most countable. For $\beta \in S^1$ let $Q_\beta = \beta+Q$. Then there exists $\beta \in S^1$ such that $Q_\beta\cap-Q_\beta = \varnothing$.
\end{lemma}
\begin{proof}
  Let $B_Q = \{\beta \in S^1\,|\,Q_\beta\cap -Q_\beta \neq \varnothing\}$. Let $\beta \in B_Q$. By definition, there are $q_\beta,q_\beta' \in Q$ such that $\beta+q_\beta = -(\beta+q_\beta')$. Fix a choice of such a pair $(q_\beta,q_\beta') \in Q^2$. This defines a map $S \fc B_Q \to Q^2$, $S(\beta) = (q_\beta,q_\beta')$. We claim that for each $(q,q')$, $S^{-1}(q,q')$ contains at most two elements. Indeed, fix $(q,q') \in Q^2$. If $S^{-1}(q,q') = \varnothing$, there is nothing to prove. Otherwise let $\beta \in B_Q$ be such that $\beta+q=-\beta+q'$, or equivalently, $2\beta=-(q+q')$. From this we see that $\beta$ is determined uniquely up to a possible addition of $\frac 12$. Thus $S^{-1}(q,q')$ equals either $\{\beta\}$ or $\{\beta,\beta+\frac12\}$. It follows that
  $$2|B_Q| \leq |\im S| \leq |Q^2|\,.$$
  Since $Q$ is at most countable, basic cardinality theory then implies that so is $B_Q$. Now any $\beta \in S^1\setminus B_Q$ will work.
\end{proof}

\begin{proof}[Proof of Theorem \ref{thm:sh_Lag_torus_contains_sh_point}]
  Let $A = \{q\in T\,|\,\mu(\{q\}) > 0\}$ be the set of atoms of $\mu|_T$. It is at most countable. Let $\sigma \in (0,\frac 12)$. By Lemma \ref{lem:decreasing_seq_meas_sets_limit}, there is $\delta \in (0,\epsilon)$ with $\mu(T\times (B_0(\delta)\setminus \{0\})) < \sigma$. Identifying $T\cong\T^n=(S^1)^n = (\R/\Z)^n$, consider the following family of functions on $T$:
  $$\cF_\delta = \left\{f(q)=\frac{\delta}{6\pi}\sum_{i=1}^{n}\sin 2\pi(q_i+\beta_i)\,\Big|\, \beta_1,\dots,\beta_n\in \R/\Z\right\}\,.$$
  For each $f \in \cF_\delta$ we have $\|df\| = \frac{\delta}{3} < \frac{\epsilon}3$. Thus $T_f = \kappa(f) + \Gamma_{df} \subset T\times \ol B_0(2\|df\|) \subset T\times B_0(\delta)$ by Remark \ref{rem:T_f_contained_ball_2_norm_df}. Moreover, $T_f\cap T$ is finite.

  We claim that there exists $f \in \cF_\delta$ such that for each $q,q' \in A$ with $q \neq q'$ we have $\pr(d_qf) \neq \pr(d_{q'}f)$. Indeed, for $i=1,\dots,n$, let $A_i\subset S^1$ be the projection of $A\subset T\cong \T^n$ onto the $i$-th circle factor. Since $A$ is at most countable, so is $A_i$, and therefore by Lemma \ref{lem:set_reflection}, there exists $\beta_i \in S^1$ such that $(\beta_i+A_i)\cap -(\beta_i+A_i) = \varnothing$. Let $f \in \cF_\delta$ be the function corresponding to a choice of such $\beta_1,\dots,\beta_n$. We claim that $f$ satisfies the required property. To show this, let $q,q'\in A$ be such that $\pr(d_qf) = \pr(d_{q'}f)$. Since $\pr\circ df$ is just the gradient of $f$, we equivalently have for each $i$: $\partial_if(q) = \partial_if(q')$, that is
  $$\cos 2\pi (q_i+\beta_i) = \cos2\pi (q_i' + \beta_i)\,.$$
  For $\theta,\theta' \in \R/\Z$ we have $\cos 2\pi\theta = \cos 2\pi\theta'$ if and only if $\theta = \theta'$ or $\theta=-\theta'$. It follows that for each $i$ we either have $q_i=q_i'$ or $q_i+\beta_i=-q_i'-\beta_i$. In the latter case we see that $q_i+\beta_i \in \beta_i + A_i$ and $-q_i'-\beta_i \in -(\beta_i+A_i)$, whence $(\beta_i+A_i)\cap - (\beta_i+A_i) \neq \varnothing$, contradicting our choice of $\beta_i$. It follows that for each $i$ we have $q_i=q_i'$, that is $q=q'$.

  Fixing such $f$, we now have
  $$\frac 12 \leq \mu(T_f) = \mu(T_f\cap T) + \mu(T_f\setminus T)\,.$$
  Since $T_f\subset T\times B_0(\delta)$, we have $T_f\setminus T\subset T\times (B_0(\delta)\setminus\{0\})$, whence $\mu(T_f\setminus T) < \sigma$. That is, we obtain
  $$\mu(T_f\cap T) > \frac 12-\sigma\,.$$
  Assume $q,q' \in A\cap T_f = A\cap(T_f\cap T)$. Since $T_f = \kappa(f) + \Gamma_{df}$, we have $q \in T_f\cap T$ if and only if $\pr(d_qf)=-\kappa(f)$, and likewise $\pr(d_{q'}f)=-\kappa(f)$. That is, $\pr(d_qf)=\pr(d_{q'}f)$, and therefore $q = q'$ by the above property of $f$. It follows that $T_f\cap T$ contains at most one element of $A$. Since $\mu(T_f\cap T) > \frac 12-\sigma >0$, $T_f\cap T$, being finite, must contain at least one point of positive measure. We conclude that $T_f\cap T$ contains exactly one element $q\in A$, and that it satisfies
  $$\mu(\{q\}) = \mu(T_f\cap T) > \frac 12 -\sigma\,.$$

  To summarize, for each $\sigma \in (0,\frac12)$ we have found $x \in A$ with $\mu(\{x\}) > \frac 12-\sigma$. Lemma \ref{lem:point_half_mass} then implies that $A$ contains an element $q$ with $\mu(\{q\}) \geq \frac12$, concluding the proof.
\end{proof}

\bibliographystyle{abbrv}
\bibliography{bibfile}

\noindent
\begin{tabular}{l}
{\bf Adi Dickstein} \\
Tel Aviv University \\
School of Mathematical Sciences\\
69978, Tel Aviv, Israel\\
\&\\
D\'epartement de Math\'ematiques et de Statistique\\
Universit\'e de Montr\'eal\\
CP 6128 succ Centre-Ville\\
Montr\'eal QC H3C 3J7\\
Canada\\
{\em E-mail:}  \texttt{adi.dickstein@gmail.com}
\end{tabular}

\medskip

\noindent
\begin{tabular}{l}
{\bf Frol Zapolsky} \\
University of Haifa \\
Department of Mathematics \\
Faculty of Natural Sciences \\
3498838, Haifa, Israel \\
\&\\
MI SANU \\
Kneza Mihaila 36 \\
Belgrade 11001\\
Serbia\\
{\em E-mail:}  \texttt{frol.zapolsky@gmail.com}
\end{tabular}

\end{document}